\documentclass[a4paper,fleqn]{cas-sc}

\usepackage{amsmath,amsthm,latexsym,amsfonts,amssymb,amsmath,amscd,epic,graphics}
\usepackage{psfrag,epsfig}
\usepackage{graphicx}

\usepackage{xcolor}

\usepackage{bm}

\newtheorem{theorem}{Theorem}[section]
\newtheorem{lemma}[theorem]{Lemma}
\newtheorem{definition}[theorem]{Definition}
\newtheorem{proposition}[theorem]{Proposition}
\newtheorem{corollary}[theorem]{Corollary}
\newtheorem{remark}[theorem]{Remark}

\newcommand{\R}{\mathbb{R}}
\newcommand{\Z}{\mathbb{Z}}
\newcommand{\Q}{\mathbb{Q}}
\newcommand{\N}{\mathbb{N}}

\newcommand{\x}{\ensuremath{\bm{x}}\xspace}
\newcommand{\y}{\ensuremath{\bm{y}}\xspace}
\newcommand{\z}{\ensuremath{\bm{z}}\xspace}
\newcommand{\w}{\ensuremath{\bm{w}}\xspace}
\newcommand{\f}{\ensuremath{\bm{f}}\xspace}
\newcommand{\g}{\ensuremath{\bm{g}}\xspace}

\renewcommand{\t}{\ensuremath{\bm{t}}\xspace}
\newcommand{\0}{\ensuremath{\bm{0}}\xspace}

\renewcommand{\b}{\ensuremath{\bm{b}}\xspace}
\renewcommand{\a}{\ensuremath{\bm{a}}\xspace}

\newcommand{\K}{\mathcal{K}}
\renewcommand{\H}{\mathcal{H}}
\renewcommand{\SS}{\mathcal{S}}

\newcommand{\dsum}{\displaystyle\sum}
\newcommand{\dmax}{\displaystyle\max}
\newcommand{\dmin}{\displaystyle\min}

\newcommand{\La}{\ensuremath{\bm{\lambda}\xspace}}

\renewcommand{\v}{\ensuremath{{\bm{v}}\xspace}}
\renewcommand{\u}{\ensuremath{{\bm{u}}\xspace}}

\newcommand{\SOCP}{\ensuremath{\rm{SOCP}}\xspace}

\newcommand{\arw}{\ensuremath{{\rm Arw}}\xspace}
\newcommand{\height}{\ensuremath{{\rm height}}\xspace}
\newcommand{\tr}{\ensuremath{\rm{Tr}}\xspace}

\def\bmx{\bm{x}}
\def\x{\bmx}
\def\a{\bm{a}}
\def\bc{\bm{c}}

\renewcommand{\;}{\quad}

\def\ip#1#2{\left\langle #1, #2 \right\rangle}

\def\norm#1{\left\|#1\right\|}
\def\norme#1{\norm{#1}_{2}}
\def\normi#1{\norm{#1}_{\infty}}

\def\lp#1#2{\norm{#1}_{#2}}

\definecolor{atomictangerine}{rgb}{1.0, 0.6, 0.4}
\definecolor{myblue}{RGB/cmyk}{0,60,101/1,.57,0,.4}

\newcommand{\rev}[1]{{\color{black} #1}}

\DeclareMathOperator{\bits}{\tau}

\usepackage{adjustbox}

\usepackage[utf8]{inputenc}  
\usepackage[T1]{fontenc}
\usepackage[numbers,sort]{natbib}
\usepackage{anyfontsize}

\begin{document}

\title[mode = title]{On the Complexity of \texorpdfstring{$p$}{}-Order Cone Programs}
\shorttitle{On the Complexity of \texorpdfstring{$p$}{}-Order Cone Programs}

\author[1]{V\'ictor Blanco}[orcid=0000-0002-7762-6461]
\ead{vblanco@ugr.es}
\author[2,3]{Victor Magron}[orcid=0000-0003-1147-3738]
\ead{vmagron@laas.fr}

\author[1]{Miguel Mart\'inez-Ant\'on}[orcid=0000-0002-4224-3655]
\ead{mmanton@ugr.es}

\affiliation[1]{organization={Institute of Mathematics (IMAG), Universidad de Granada}, country={Spain}}
\affiliation[2]{organization={LAAS-CNRS}, country={France}}
\affiliation[3]{organization={Institute of Mathematics of Toulouse}, country={France}}

\shortauthors{Blanco, Magron, and Mart\'inez-Ant\'on}

\begin{abstract}
This manuscript explores novel complexity results for the feasibility problem over $p$-order cones, extending the foundational work of Porkolab and Khachiyan [J. Glob. Optim. 10(4) (1997), pp. 351--365]. 
By leveraging the intrinsic structure of $p$-order cones, we derive refined complexity bounds that surpass those obtained via standard semidefinite programming reformulations. 
Our analysis not only improves theoretical bounds but also provides practical insights into the computational efficiency of solving such problems. 
In addition to establishing complexity results, we derive explicit bounds for solutions when the feasibility problem admits one. For infeasible instances, we analyze their discrepancy quantifying the degree of infeasibility. 
Finally, we examine specific cases of interest, highlighting scenarios where the geometry of $p$-order cones or problem structure yields further computational simplifications. 
These findings contribute to both the theoretical understanding and practical tractability of optimization problems involving $p$-order cones.
\end{abstract}


\begin{keywords}
$p$-order cone, semidefinite programming, computational complexity, discrepancy.
\end{keywords}

\makeatletter\def\Hy@Warning#1{}\makeatother
\maketitle

\section{Introduction}
Second Order Cones (SOC) are fundamental structures in conic optimization, that have been proven to be useful in different applied disciplines such as engineering design~\citep{blanco2024optimal,goez2019second,mohammadisiahroudi2024generating}, robust optimization~\citep{zhen2022robust},  financial modeling~\citep{brar2020portfolio,lu2006new}, or machine learning~\citep{lopez2016multi,kucukyavuz2023consistent,maldonado2014imbalanced}, among many others~\citep[see e.g.][]{lobo1998applications,alizadeh2003second}. 
Recently a SOC representation of a specific cone of sparse nonnegative polynomials has been proposed in  \cite{averkov2019optimal} and \cite{magron2023sonc}. 
The SOC (also known as the Lorentz cone) in $\R^{n+1}$ is defined as $\K_2^n = \{(\x,z) \in \R^{n+1}: \|\x\|_2\leq z\}$. The natural extension of these cones are the $p$-order cones, where the Euclidean norm is replaced by the $\ell_p$-norm. These cones have been recommended in different applications where measuring with \rev{non-Euclidean}-based distances is recommended. For instance, in \citep{blanco2020lp}, Support Vector Machines with $\ell_p$-norm based margins are analyzed. In their computational experiments, the authors conclude that some datasets obtain the benefit of these norms with respect to the classical (euclidean) classifiers in terms of the obtained accuracies and the number of non-zero features used for the classification. In multiple criteria and combinatorial optimization, $\ell_p$-based aggregation of the objectives has resulted to be beneficial in terms of fairness~\citep{bektacs2020using,blanco2023fairness,kostreva2004equitable}. The incorporation of $p$-order cones into any of these optimization problems is performed, for practical purposes by \emph{rewriting} these cones as a finite set of SOCs~\citep{alizadeh2003second,blanco2024minimal}, which are the cones that the available optimization solvers allow. 

In this paper, we analyze crucial complexity questions that arise on $p$-order cones optimization, which, as far as we know have not been previously studied. Specifically, we explore the theoretical complexity of the \emph{$p$-order feasibility problem}. Thus, the main object under study is the $(n+1)$-dimensional $p$-order cone,  which is defined as:
$$
    \K^n_p:=\{(\x,z) \in \R^{n+1} : \lp{\x}{p}\leq z\},
$$
for any rational $p=\frac{r}{s}$, with $r > s\in \N^*$ and $\gcd(r,s)=1$,  and where $\|\cdot\|_p$ stands for the $\ell_p$-\emph{norm} on $\R^n$ \rev{($\|\x\|_p = \left(\sum_{i=1}^n |x_i|^p\right)^{\frac{1}{p}}$, for $\x \in \R^n$)}. 
The $p$-order feasibility  problem {\rev{(\ref{eq:OCFP}, for short)}} that consists of certifying the existence of a solution of a linear system where the variables belong to $\K_p^n$. For the case $p=2$, this convex set can be seen as a subset of the cone of $(n+1)\times (n+1)$ real symmetric positive semidefinite matrices,  $\SS^{n+1}_+$.  Thus, the feasibility problem for this case is a particular case of the \emph{semidefinite feasibility problem} already studied in  \cite{porkolab1997complexity}, by direct identification using the Schur complement. Hence, the case $p=2$ of \ref{eq:OCFP} can be analyzed using the results in ~\citep{porkolab1997complexity} for the semidefinite feasibility problem.

In this paper, we derive several complexity results for the $p$-order feasibility problem, measured in terms of the number of required arithmetic operations and the bit-length of the numbers involved. \rev{We consider the setting where the problem includes $m$ linear constraints with $n+1$ variables and coefficients of bit size at most $\tau$, along with a single $p$-order cone constraint.} In Table \ref{table:soa2} we summarize the main results derived in this paper, indicating the numbered Lemma or Theorem where the result is stated and proved. In the first part of the paper we derive the complexity results by \emph{rewriting} the problem equivalently as a SDP feasibility problem, and then, applying the complexity result in \cite{porkolab1997complexity} (Theorem \ref{lem:sdp}). Corollary \ref{lem:int-points} is a consequence of using this result together with the interior-point complexity analysis provided in \cite{ben2001lectures} for the case when the coefficients of the system have at most bit size $\tau$. Nevertheless, in this paper we  explicitly exploit the structure of the cone $\K_p^n$, avoiding the rewriting of this cone as a subset of the semidefinite cone, but still following the same paradigm of \cite{porkolab1997complexity}. Thus, we get improved complexity bounds for the general case (Theorem \ref{th:complexgc}), but, additionally, we analyze particular cases for the value of $p$, where the complexity can be further improved, as the case when $p$ is even (Theorem \ref{th:complexeven}), when $p$ is in the form $\frac{r}{r-1}$ with $r$ even (Corollary \ref{cor:complexdualeven}), or the very special case of the SOCP (Theorem \ref{th:socpcomplexity}).
\begin{table}[H]
\centering
\begin{tabular}{>{\centering\arraybackslash}m{4cm} >{\centering\arraybackslash}m{2.5cm} >{\centering\arraybackslash}m{6cm}}
\hline
 & $p$ & No. arithmetic operations \\
\hline\hline
Theorem \ref{lem:sdp} & $\frac{r}{s}\in \Q_{> 1}$ & $m(n\log r)^{O((n\log r)^2)}$ \\
\hline
Corollary \ref{lem:int-points} & $2$ (SOCP) & $\bits m n^3 \min\{m,n\}^{O(\min\{m,n\})}$ \\
\hline\hline
Theorem \ref{th:complexgc}  & $\frac{r}{s}\in \Q_{> 1}$ & $\min\{[r(m+n)]^{O(n)}, m(rn)^{O(n^2)}\}$ \\
\hline
Theorem \ref{th:complexeven}& $r\in 2\N^*$ & $m(r\min\{m,n\})^{O(n)}$ \\
\hline
Corollary \ref{cor:complexdualeven} & $\frac{r}{r-1}, r\in 2\N^*$ & $m(r\min\{m,n\})^{O(\min\{m,n^2\})}$ \\
\hline
Theorem \ref{th:socpcomplexity} & $2$ (SOCP) & $m\min\{m,n\}^{O(\min\{m,n\})}$ \\
\hline
\end{tabular}
\caption{Summary of our complexity results for the $p$-order feasibility problem}\label{table:soa2}
\end{table}

Additionally, we use the geometrical and algebraic structure of $p$-order cone to derive bounds on either the $\log$-modulus of feasible solutions, or the discrepancy otherwise, and bounds on the number of arithmetic operations and bit size of the involved numbers for the feasibility problem. In case the $p$-order cone problem is feasible we derive bounds for Euclidean norm of a solution $(\x,z)$ to the system. In case the problem is infeasible, we compute bounds for the minimal violation of the linear system involved in the feasibility problem. Table \ref{table:disc} summarizes the results we obtained in this line. In the first part we detail the logarithmic upper bound for a feasible solution to the problem. In the second part, we detail the  logarithm of the bound for the minimal violation (the so-called $\log$-discrepancy).

\begin{table}[h!]
\centering
\begin{tabular}{>{\centering\arraybackslash}m{4cm} >{\centering\arraybackslash}m{2.5cm} >{\centering\arraybackslash}m{6cm}}
\hline
 & $p$ & $\log$ norm\\
\hline\hline
Lemma \ref{th:boundseven} & $r \in 2\N^*$ & $\bits(\min\{m,n\}r)^{O(n)}$  \\
Theorem \ref{th:bounds} & $\frac{r}{s} \in \Q_{>1}$ & $\bits(nr)^{O(n\min\{m,n\})}$ \\
Theorem \ref{th:bounds2} & $\frac{r}{r-1}$, $r\in 2\N^*$ & $\bits(r\min\{m,n\})^{O(\min\{m,n\})}$\\
\hline
& $p$ & $\log$ discrepancy \\\hline\hline
Theorem \ref{th:discrepancy} &  $\frac{r}{s} \in \Q_{>1}$ & $-\bits(rn)^{O(n\min\{m,n\})}$\\
Theorem \ref{th:discrepancy} &  $\frac{r}{r-1}$, $r\in 2\N^*$ & $-\bits(r\min\{m,n\})^{O(\min\{m,n\})}$\\\hline
\end{tabular}
\caption{Summary of the bounds on feasible solutions and discrepancies for $p$-order cone feasibility problems.}\label{table:disc}
\end{table}

The complexity results obtained in this paper have a direct impact in the complexity of solving a $p$-order cone optimization problems, i.e., minimizing a linear function subject to some linear constraints and the requirements that some of the variables belong to $\K_p^n$. Using the feasibility algorithm as an oracle, one can apply any binary search technique to derive, up to some tolerance $\epsilon>0$, the optimal value for the optimization problems. Some complexity results in this line will be derived for some interesting families of optimization problems involving $\ell_p$-norms, as the norm minimization problem, support vector machines, continuous location problems, robust least squares, or robust linear optimization.

The rest of the paper is organized as follows. In Section \ref{sec_prelim} we state the notation and main results \rev{used} to derive the complexity results in this paper. In this section we also prove the complexity results that can be derived for the problem by rewriting the problem as a semidefinite programming problem. In Section \ref{sec:disc} we provide the upper bounds for the feasible solutions and the discrepancy for the $p$-order feasibility problem. Section \ref{sec:comp} is devoted to prove the complexity bounds for the $p$-order feasibility problem. Exploiting the geometry of $p$-order cones we apply a similar strategy as the one by Porkolab and Khachiyan to derive new complexity bounds for the problem. 
We also analyze some particular cases for the value of $p$, namely, $p=2$, $p$ nonnegative even integer, and $p$ of the form $\frac{r}{r-1}$ with $r$ nonnegative integer. 
In Section \ref{sec:multi}, we extend our complexity results to problems involving several $p$-order cone constraints. 
Finally, in Section \ref{sec:applis} we emphasize how the complexity bounds are applied to optimization problems that involve $p$-order cones.

\section{Preliminaries}
\label{sec_prelim}
In this section we set the notation for the rest of the sections and state the main results that will be useful in our developments.

\paragraph*{\bf Notation for complexity estimates}  
For complexity estimates, we use the bit complexity model. 
For an integer $b \in \Z \backslash \{0\}$, we denote by $\tau(b) := \lfloor \log_2 (|b|) \rfloor +1$ the bit size of $b$, with the convention $\tau(0) := 1$. 
For a rational number $a = \frac{b}{c}$, with $b \in \Z, c \in \Z \backslash \{0\}$ and ${\rm gcd}(b,c)=1$, we denote $\max \{\tau (b), \tau (c)\}$ by $\tau (a)$. 
For the sake of simplicity, all derived complexity estimates are provided while assuming integer input data. It is straightforward to establish similar estimates while assuming rational input data. 

For two mappings $g, h : \N^l \to \R_{>0}$, the expression ``$g(v) = O(h(v))$'' means that there exist integers $b,N \in \N$ such that when all coordinates of $v$ are greater than or equal to $N$, $g(v) \leq b h(v)$. 

As already mentioned, the main question that we address in this paper is on the complexity of the $p$-order feasibility problem which is defined as follows.
\begin{definition}[$p$-order feasibility problem]
The \emph{$p$-order feasibility problem} consists of determining whether there exists a real vector $(\x,z)\in\R^{n+1}$ such that
 \begin{align}
     F \x + G z \leq H, &  \nonumber\\
     (\x,z) \in \K_p^n. &\label{eq:OCFP}\tag{${\bf F}_p$}
 \end{align}
 where $F \in \Z^{m\times n}, G, H \in \Z^m$ for $m\in \N^*$.
\end{definition}

We will denote by $\f_i^T$ the $i$th row of matrix $F$, for  $i=1,\ldots,m$, and  $G=(g_1,\ldots, g_m)$, and $H=(h_1,\ldots, h_m)$. 

Note that when $n=1$ the $\ell_p$-norm becomes $|\cdot|$ and the feasibility problem \ref{eq:OCFP} is trivial. Therefore we assume throughout the paper that $n\geq 2$.

The special case of the SOCP ($p=2$) can be seen as a particular case of the \emph{semidefinite feasibility problem} by a direct identification of the Lorentz cone with a subcone of the SDP cone using the Schur complement. Thus, the complexity of the $2$-order feasibility problem can be analyzed using the results in \citep{porkolab1997complexity} for the semidefinite feasibility problem which is defined by the existence of a real $n\times n$ symmetric matrix $X\in \SS^n_+$ which is solution for
\begin{equation}\label{eq:SFP}
   \ip{A_i}{X} \leq b_i, \; i=1,\ldots,m; \; X\succeq 0, \tag{${\bf F}_{\rm SDP}$}
\end{equation}
where $A_1,\ldots, A_m\in \SS^n$ are integral $n\times n$ symmetric matrices, $b_1,\ldots, b_m\in\Z$, $\ip{A}{X}:=\tr(AX)$ denotes the standard inner product on the space of real symmetric matrices and $X\succeq 0$ stands for the membership of the symmetric matrix $X\in\SS^n$ in the cone of positive semidefinite matrices. 

On the other hand, \cite{blanco2024minimal} derived minimal SOC reformulations of the $p$-order cone. Thus, combining these results and those for the SDP \rev{feasibility} problem, the complexity results obtained for the SOCP can be extended to any $p\in \Q_{> 1}$ by applying an explicit and minimal semidefinite extended representation of \ref{eq:OCFP}, and then use of the results in ~\citep{porkolab1997complexity}, as we prove in the following result.
\begin{theorem}
    \label{lem:sdp}
Let $p = \frac{r}{s} \in \Q$ with $r > s \in \N^*$ and $\gcd(r, s) = 1$. Then, the feasibility of~\ref{eq:OCFP} can be tested using $m(n \log r)^{O((n \log r)^2)}$ arithmetic operations over numbers of bit size $\bits (n \log r)^{O((n \log r)^2)}$\rev{, where $\bits$ denotes the maximum bit size of the coefficients.}
\end{theorem}

\begin{proof}
Let us first analyze the case $p=2$, and then we extend the result for any $p$. For each vector $(\x,z)\in\R^{n+1}$, the \emph{arrow-shaped} matrix $\arw (\x,z)\in\SS^{n+1}$ is defined as:
\begin{equation*}
    \arw(\x,z):=\begin{bmatrix}
    zI_n & \x\\
    \x^T & z
\end{bmatrix}.
\end{equation*}
Note that by the Schur complement, $(\x,z)\in \K^n_2$ if and only if $\arw(\x,z)\succeq 0$. Thus, \ref{eq:OCFP}, with $p=2$, is equivalent to
{\small
\begin{equation}
    \ip{\arw\left(\frac{1}{2}\f_i,\frac{g_i}{n+1}\right)}{\arw(\x,z)} \leq h_i \; i=1,\ldots,m, \; \arw(\x,z)\succeq 0.\label{eq:equiv2SDP}
\end{equation}}
Based on the identification of ${\bf F}_2$ as \eqref{eq:equiv2SDP}, to rewrite this feasibility problem as \ref{eq:SFP} one must note that in the condition $\arw(\x,z)\succeq 0$ are hidden $\frac{1}{2}n(n+1)$ extra linear constraints. Hence, by ~\cite{porkolab1997complexity}, if the integer coefficients of \ref{eq:OCFP} have at most bit size $\tau$, its feasibility can be tested in $mn^{O(n^2)}$ arithmetic operations over $\bits n^{O(n^2)}$-bit numbers. Continuing with the general case. We also can define the three-dimensional rotated second order cone as
\begin{equation*}
    \K^2\left(\frac{1}{2}\right):=\{(x,y,z)\in \R^3 : x^2\leq yz, \; y,z\geq 0\}.
\end{equation*}
This cone is obtained from $\K^2_2$ after a rotation, note that $(x,y,z)\in \K^2\left(\frac{1}{2}\right)$ if and only if $(2x,y-z,y+z)\in \K^2_2$. As already mentioned $\K^2_2$ is a subset of \rev{the cone $\SS^3_+$ of $3\times 3$ real positive semidefinite matrices}. But, if we consider the three-dimensional rotated second order cone it is easy to see that $\K^2\left(\frac{1}{2}\right)=\SS^2_+$ \rev{(the cone of $2\times 2$ real positive semidefinite matrices)}. 

 In \cite[][Theorem 15]{blanco2024minimal}, the authors derived a minimal three-dimensional rotated second order cone ($\SS^2_+$) extended representation based on a minimal (in the sense of number of vertices) graph structure named mediated graph. In the case of the $p$-order cone with $p=\frac{r}{s}$, $r > s$, and $\gcd(r,s)=1$ this graph $G=(V\cup\{0,r\},A)$ satisfies $V\subset (0,r)$, $|V|=\lceil\log r\rceil$, $s\in V$, and it is defined by the sets of outgoing arcs $\delta^+(v)=\{u,2v-u\}\subset V\cup\{0,r\}$ for all $v\in V$. With all of this, \ref{eq:OCFP} is equivalent to
 \begin{align*}
    & F^T\x + G^Tz \leq H, & \\
    & \dsum_{j=1}^{d} t_ j -z\leq 0, & \\
    & {\rm BlockDiag}\left( W_{jv} \right)_{\substack{j=1,\ldots,n \\ v \in V}}
 \succeq 0, &
\end{align*}
where $W_{jv}=\begin{bmatrix}
    w_{ju} & w_{jv}\\
    w_{jv} & w_{j,2v-u}
\end{bmatrix}\in \SS^2$, with $j=1,\ldots,n$; $v\in V$, and $u, 2v-u\in \delta^+(v)$; and $w_{j0}=z$, $w_{js}=x_j$, and $w_{jr}=t_j.$ Again, in the semidefinite constraint are hidden $O((n\log r)^2)$ affine constraints thus, by ~\cite{porkolab1997complexity}, its feasibility can be tested in $m(n\log r)^{O((n\log r)^2)}$ arithmetic operations over $\bits (n\log r)^{O((n\log r)^2)}$-bit numbers.
\end{proof}
Besides, for $p=2$ (SOCP), we can use our results combined with the complexity of the interior-points methods where the number of arithmetic operations to solve ${\bf F}_2$ problem with $n$ variables and $m$ constraints at accuracy $\varepsilon$ stated in \cite[\S~4.6.2]{ben2001lectures} is $ O\left(\sqrt{(m+1)} n (n^2 + m) \log(1/\varepsilon) \right)$.  
\begin{corollary}
    \label{lem:int-points}
     $({\bf{F}_2})$ can be tested in $\bits m n^3 \min\{m,n\}^{O(\min\{m,n\})}$ arithmetic operations using interior-points methods.
\end{corollary}
\begin{proof}
    Plugging $\log(1/\varepsilon) = \bits \min\{m,n\}^{O(\min\{m,n\})}$ (see Theorem \ref{th:bounds2}) provides a bound on how many operations are required to test feasibility that after reductions become in $\bits m n^3 \min\{m,n\}^{O(\min\{m,n\})}$ where $\bits$ is the maximum of the bit size in the input of the problem.
\end{proof}
In the previous results, the complexity of \ref{eq:OCFP} is obtained by rewriting the problem as a particular case of \ref{eq:SFP} and then apply the results for the general case. As already announced, better complexity bounds can be obtained by explicitly exploiting the structure of $\K_p^n$. The following definitions and results will allow us to derive these new complexity results.
\begin{definition}[Dual Cone]
    Let $\K\subset \R^n$ be a convex cone. Its dual cone is defined as
$$
    \K^*:=\{\y\in \R^n : \ip{\y}{\x} \geq 0 \; \forall \x\in \K\}.
$$
\end{definition}
The dual cone of the $p$- order cone $\K^n_p$ is the $q$-order cone $\K^n_q$, where $q = \frac{p}{p-1}$ (the \emph{conjugate} of $p$,  i.e. $\frac{1}{p}+\frac{1}{q}=1$). For $p=q=2$, one has $(\K^n_p)^* = \K^n_q = \K^n_p$, so   the second order cone is said self-dual. 
\begin{definition}
  Let $A$ be a nonempty set in a vector space $X$. The recession cone of $A$ is defined as: 
$$
    {\rm recc}(A):=\{\y\in X:\forall \x\in A,\forall \lambda \geq 0, \x+\lambda \y\in A\}.
$$
\end{definition}
Note that if $A$ is the intersection of a proper cone and an affine space as \ref{eq:OCFP}, then the recession cone is equivalent to

\begin{equation}\label{eq:reccbis}
        \K^n_p \cap \bigcap_{i=1}^m \H_i,
    \end{equation}
where $\H_i$ is the halfspace $\{(\x,z)\in \R^{n+1} : \f^T_i\x+g_iz \leq 0\}$. If $A$ is a nonempty closed convex subset of a finite-dimensional Hausdorff space (e.g. $\R^n$), then ${\rm recc}( A ) = \{ \0 \}$ if and only if $A$ is bounded.

The main tools that we apply to derive the complexity results for \ref{eq:OCFP} come from first-order theory of the reals and the results in \cite{renegar1992computationalb}. In what follows, we recall the fundamental results required in our proofs.
\begin{definition}
Given $k,\omega, d \in \N$, a formula in the first-order theory of the reals is an expression of the form
\begin{equation}\label{eq:sf}
    (Q_1\x_1\in \R^{n_1})\ldots (Q_\omega \x_\omega \in \R^{n_\omega})P(\y,\x_1,\ldots,\x_\omega), \tag{SF}
\end{equation}
where: 
\begin{itemize}
    \item $\y\in \R^k$ is a free variables vector;

    \item each $Q_l$ ($l=1,\dots,\omega$) is one of the quantifiers $\forall$ or $\exists$;

    \item $P(\y,\x_1,\ldots,\x_\omega)$ is a quantifier free boolean formula with $m$ atomic predicates of the form
$$
        g_j(\y,\x_1,\ldots, \x_\omega)\lhd_j 0, \; j=1,\ldots,m,
$$
    where $\lhd_j \in \{<,>,\leq,\geq,=,\neq\}$ and $g_j$ is a real polynomial of degree at most $d$.
\end{itemize}
\end{definition}
Note that the above formula is in prenex form, i.e., all quantifiers in \ref{eq:sf} appear in front. Formulas with no free variables are called sentences. We say $\y\in\R^k$ is a solution of \ref{eq:sf} if the sentence obtained by substituting $\y$ into \ref{eq:sf} is true.
\begin{proposition}[\cite{renegar1992computationalb}, Proposition 1.3]\label{prop:renegarbound}
    If a formula \ref{eq:sf} has only integer coefficients, each of bit size at most $\bits$, then every connected component of the set of its solutions intersects the ball $\{\y\in \R^k : \norme{\y}\leq R\}$, where $R$ satisfies $
        \log R\leq \bits(md)^{2^{O(\omega)}k\prod_{l=1}^\omega n_l}$.
\end{proposition}
\begin{theorem}[\cite{renegar1992computational}, Theorem 1.2]\label{prop:renegarfree}
    There is an algorithm which, given a formula \ref{eq:sf}, finds an equivalent quantifier free formula of the form
$$
        \bigvee_{i=1}^I \bigwedge_{j=1}^{J_i}(h_{ij}(\y)\lhd_{ij} 0 ),
$$
    where:
    \begin{itemize}
        \item $I\leq (md)^{2^{O(\omega)}k\prod_{l=1}^\omega n_l}$,
        \item $J_i\leq (md)^{2^{O(\omega)}k\prod_{l=1}^\omega n_l}$,
        \item  $\deg h_{ij}(\y)\leq (md)^{2^{O(\omega)}k\prod_{l=1}^\omega n_l}$.
    \end{itemize}
    The algorithm requires $(rd)^{2^{O(\omega)}k\prod_{l=1}^\omega n_l}$ operations and $(rd)^{O(k+\sum_{l=1}^\omega n_l)}$ evaluations of the input formula. 
    If the coefficients of the atomic polynomials $g_j$, $j=1,\ldots, r$, are integers of bit size at most $B$, then the algorithm works with numbers of bit size $
        (B+k)(md)^{2^{O(\omega)}k\prod_{l=1}^\omega n_l}$. This bound also holds for the bit size of the coefficients of polynomials $h_{ij}$.
\end{theorem}
\begin{theorem}[\cite{renegar1992computational}, Theorem 1.1]\label{prop:renegarcomplexity}
    There is an algorithm for the decision problem of the first-order theory of the reals that requires $(md)^{2^{O(\omega)}k\prod_{l=1}^\omega n_l}$ 
    operations and $(md)^{O(\sum_{l=1}^\omega n_l)}$ evaluations of the input formula. When restricted to sentences involving only polynomials with integer coefficients of bit size at most $B$, the procedure works with numbers of bit size $B(md)^{2^{O(\omega)}k\prod_l n_l}$.
\end{theorem}
The following inequality is a well-known bound on nonzero roots of univariate polynomials~\citep[see, e.g.][p.~261]{mignotte1982some}.
\begin{proposition}\label{prop:height}
    Let $p(x)=\dsum_{k=0}^da_kx^k$ be a univariate polynomial with integer coefficients, and $\alpha$ be a nonzero root of $p$. Then $|\alpha|\geq \frac{1}{1+h}$, where $h=\dmax_{0\leq k\leq d}\{|a_k|\}$ is the height of $p$.
\end{proposition}
The result below is a suitable variant of the Fundamental Theorem of Linear Inequalities ~\citep[see, e.g.][]{schrijver1998theory}.
\begin{proposition}\label{prop:FTSLI}
    Consider a system of linear inequalities:
    $$
        \a_i^T\x\leq b_i, \; i\in M=\{1,\ldots,m\},
    $$
    and let $K\subset \R^n$ be a convex set. If $P=K\cap \{\x\in\R^n : \a_i^T\x\leq b_i, i\in M\}$ is nonempty, then there exists a subset $I\subseteq M$ such that $|I|\leq \min\{m,n\}$ and $\emptyset\neq K\cap \{\x\in\R^n : \a_i^T\x\leq b_i, i\in I\}\subseteq P$.
\end{proposition}

\section{Bounding Solutions and  Discrepancies}\label{sec:disc}
In this section we derive some results related to the known answer to \ref{eq:OCFP}. Specifically, in case the problem is certified to be feasible we derive upper bounds on the log-modulus of at least one feasible solution. In case the problem is infeasible, we analyze its discrepancy, that is the minimum violation of the feasibility. We also derive lower bounds to this value in case the system is infeasible.

Given a positive $m \in \N$, we denote the $(m-1)$-dimensional simplex as:
\begin{equation}\label{def:simplex}
    \Delta_m := \left\{\La \in \R^m : \dsum_{i=1}^m \lambda_i =1, \lambda_1\geq0,\ldots, \lambda_m\geq 0\right\}.
\end{equation}
The first result that we address is the one of representing $\ell_p$-norm based constraints in the first order formula language.
\begin{lemma}\label{th:exreprs}
Let $p=\frac{r}{s}$  with $r > s \in \N^*$ and $\gcd(r,s)=1$. 
   \rev{ \begin{enumerate}
        \item $\lp{\x}{p}= z$ if and only if $(\x, z)$ satisfies the following first-order formula:
        \begin{multline}\label{eq_gpcf_1}
        \exists \t \in \R^{n} \left\{\bigwedge_{j=1}^n \left[ (x_j^r= z^{r-s}t_j^s) \vee (-x_j^r= z^{r-s}t_j^s)\right], 
        \dsum_{j=1}^n t_j = z, t_1\geq 0,\ldots, t_n\geq 0\right\}.
        \end{multline}
        \item $\lp{\x}{p}\leq z$ if and only if $(\x, z)$ satisfies the following first-order formula:
        \begin{multline}\label{eq_gpcf_2}
        \exists \t \in \R^{n} \left\{\bigwedge_{j=1}^n \left[ (x_j^r\leq z^{r-s}t_j^s) \wedge (-x_j^r\leq z^{r-s}t_j^s)\right], 
        \dsum_{j=1}^n t_j \leq z, t_1\geq 0,\ldots, t_n\geq 0\right\}.
    \end{multline}
    \end{enumerate}}
    
\end{lemma}
\begin{proof}
    The proof follows by rewriting equivalently the inequality/equation $\lp{\x}{p}\lhd z$, \rev{where $\lhd \in \{\leq, =\}$}, as a standard formula \ref{eq:sf} as detailed in ~\citep{blanco2024minimal}.
\end{proof}
\begin{lemma}
\label{lemma:ptheq}
Given a positive $p \in \Q$ with conjugate $q$ one has
\begin{align}
    & \dmin_{(\x,z)\in \K^n_p}\{\f^T\x + gz : z=R\}=R(g-\lp{\f}{q}), \label{eq:min}\\
    & \dmin_{(\x,z)\in \K^n_p}\{\f^T\x + gz : z\leq R\}=\min\{0,R(g-\lp{\f}{q})\}.  \label{eq:min2}
\end{align}
\end{lemma}
\begin{proof}
To show the first identity, observe that
\begin{align*}
    \dmin_{(\x,z)\in \K^n_p}\{\f^T\x + gz : z=R\}
        &=R(g-\lp{\f}{q}) + \dmin_{(\x,z)\in \K^n_p}\{\f^T\x + \lp{\f}{q}z : z=R\}
         = R(g-\lp{\f}{q}),
\end{align*}
where the last equality follows from the fact that $(\f,\lp{\f}{q})$ belongs to the boundary of $\K^n_q=(\K^n_p)^*$. \\
For the second identity, note that if $(\f,g)\in \K^n_q$, then $\dmin_{(\x,z)\in \K^n_p}\{\f^T\x + gz : z\leq R\}=0$. Otherwise $g-\lp{\f}{q}<0$, which means that the minimum on the left-hand side of \eqref{eq:min2} is negative and hence it is attained at a vector $(\x,R)$. Then \eqref{eq:min2} becomes a consequence of \eqref{eq:min}.
\end{proof}
In what follows we analyze the bounds on the feasible solutions of \ref{eq:OCFP}. We first analyze the case  $p=r\in 2\N^*$. 
\begin{lemma}\label{th:boundseven} 
Let $p=r\in 2\N^*$. Then, there exists a feasible solution $(\x,z)$ of \ref{eq:OCFP}  such that $\lp{(\x,z)}{2}\leq R$, where $\log R=\bits(\min\{m,n\}r)^{O(n)}$.
\end{lemma}
\begin{proof}
In this case the problem \ref{eq:OCFP} is equivalent to the quantifier free formula:
\begin{equation}
\label{eq:evenrSF}
    \{\f_1^T\x+g_1z-h_1\leq 0,\ldots, \f^T\x_m+g_mz-h_m\leq 0, x_1^r+\cdots +x_n^r\leq z^r, z\geq 0\}.
\end{equation}
Observe that this formula consists of $m+2$ polynomial inequalities of degree at most $r$ in $(n+1)$ free variables. If the integer coefficients appearing in \eqref{eq:evenrSF} have height at most $2^{\bits}$, then Proposition \ref{prop:renegarbound} implies that any feasible solution satisfies $\norme{(\x,z)}\leq R$, where $\log R\leq \bits(mr)^{O(n)}$.

By Proposition \ref{prop:FTSLI}, there is a set $I\subseteq M$ of size at most $\min\{m,n+1\}$ such that the system
    \begin{equation*}
        \f^T_i\x+g_iz=h_i, \; i\in I, \; \lp{\x}{p}\leq z,
    \end{equation*}

    is feasible, and any of its solutions solve the original problem \ref{eq:OCFP}. For this reason, we can obtain a better bound by replacing $m$ with $\min\{m,n+1\}$, namely $\log R\leq \bits(\min\{m,n\}r)^{O(n)}$.
\end{proof}
The next result that we prove is the general case, which is the analog of~\cite[Theorem 3.1]{porkolab1997complexity} for \ref{eq:SFP}.
\begin{theorem}\label{th:bounds} 
Let $p=\frac{r}{s} \in \Q$, with $r > s\in \N^*$ and $\gcd(r,s)=1$. 
Then:
    \begin{enumerate}
        \item There exists a feasible solution $(\x,z)$ of \ref{eq:OCFP}  such that $\lp{(\x,z)}{2}\leq R$, where $\log R=\bits(nr)^{O(n\min\{m,n\})}$. \label{th:boundpart1}

        \item Moreover, if the feasible set of \ref{eq:OCFP} is bounded, then the above bound holds for any feasible solution of \ref{eq:OCFP}. \label{th:boundpart2}
    \end{enumerate}
\end{theorem}
\begin{proof}
    Suppose that the problem \ref{eq:OCFP} is feasible, and let us define:
    \begin{align*}
        \Omega_R& :=\{(\x,z)\in \K^n_p : z=R\}, \\
        \Theta(R)& := \dmin_{(\x,z)\in \Omega_R}\max\{\f_1^T\x + g_1 z -h_1,\ldots, \f_m^T\x + g_mz - h_m\}.
    \end{align*}
    For any $R\geq 0$ and considering $\Delta_m$ as in \eqref{def:simplex},
    \begin{multline*}
        \Theta(R)=\dmin_{(\x,z)\in \Omega_R}\dmax_{\La\in \Delta_m} \dsum_{i=1}^m \lambda_i(\f_i^T\x + g_iz - h_i) \\
        =\dmax_{\La\in \Delta_m}\dmin_{(\x,z)\in \Omega_R} \left[\left(\dsum_{i=1}^m\lambda_i\f_i^T\right)\x + \left(\dsum_{i=1}^m \lambda_i g_i\right)z - \dsum_{i=1}^m \lambda_i h_i\right] 
        =\dmax_{\La\in \Delta_m}\left[ R\left(\dsum_{i=1}^m \lambda_i g_i - \lp{\dsum_{i=1}^m\lambda_i\f_i}{q}\right)  - \dsum_{i=1}^m \lambda_i h_i \right].
    \end{multline*}

    The second equality follows from Von Neumann's minimax theorem (see, e.g.~\cite{rockafellar1970convex}) and the last one follows the first identity \eqref{eq:min} from Lemma \ref{lemma:ptheq}. 
    
    Now, consider the formula 
    \begin{equation*}
        \Phi(R) := \forall \La \in \Delta_m \left\{ R\left(\dsum_{i=1}^m \lambda_i g_i - \lp{\dsum_{i=1}^m\lambda_i\f_i}{q}\right)  - \dsum_{i=1}^m \lambda_i h_i\leq 0 \wedge R\geq 0\right\}.
    \end{equation*}
    By Lemma \ref{th:exreprs} \eqref{eq_gpcf_1}, $\Phi(R)$ can be rewritten in the standard form \ref{eq:sf} as follows
    \begin{multline*}
    \forall \La \in \R^m \; \exists (\t,w) \in \R^{n+1} \Bigg\{ \Bigg\{ \left[
    \lambda_1 \geq 0, \ldots, \lambda_m \geq 0, \quad \sum_{i=1}^m \lambda_i = 1
    \right] \implies \\
    \Bigg[
    \sum_{j=1}^n t_j = w, \bigwedge_{j=1}^n \left[\left[
        \left( \sum_{i=1}^m \lambda_i f_{ij} \right)^r = w^s t_j^{r-s}
    \right]\vee \left[
        -\left( \sum_{i=1}^m \lambda_i f_{ij} \right)^r = w^s t_j^{r-s}
    \right]\right], \\
    R \left( \sum_{i=1}^m \lambda_i g_i - w \right) - \sum_{i=1}^m \lambda_i h_i \leq 0, t_1\geq 0,\ldots, t_n\geq0
    \Bigg]\Bigg\} \wedge (R \geq 0)\Bigg\}.
    \end{multline*}
Then, for any $R\in\R$, the following statements are equivalent:
    \begin{itemize}
        \item \ref{eq:OCFP} has a feasible solution in $\Omega_R$,

        \item $\Theta(R)\leq 0$,

        \item $R$ satisfies $\Phi(R)$.
    \end{itemize}
    By our original assumption, \ref{eq:OCFP} is feasible, and hence there is a nonnegative $R$ that satisfies $\Phi(R)$. 
    Next, $\Phi(R)$ is a standard formula \ref{eq:sf} of degree at most $r$ with $k=1$ free variable and $\omega=2$ quantifiers. 
    Furthermore, $\Phi(R)$ consists of $O(m+n)$ atomic \rev{polynomial} inequalities in $O(m+n)$ variables. 
    The expression $\left( \sum_{i=1}^m \lambda_i f_{ij} \right)^r$ involves integer coefficients of height at most $2^{\bits}$, thus the expansion yields coefficients with bit size at most $B=r(\bits+\log(m))$.

    Now from Proposition \ref{prop:renegarbound} it follows that $\Phi(R)$ can be satisfied by a positive number $R$ such that
    \begin{equation}\label{eq:bound}
        \log R = r(\bits+\log(m))(O(m+n)r)^{O(mn)}=\bits((m+n)r)^{O(mn)}.
    \end{equation}
    By Proposition \ref{prop:FTSLI}, there is a set $I\subseteq M$ of size at most $\min\{m,n+1\}$ such that the system
    \begin{equation*}
        \f^T_i\x+g_iz=h_i, \; i\in I, \; \lp{\x}{p}\leq z,
    \end{equation*}
    is feasible, and any solution solves the original problem~\ref{eq:OCFP}. 
    For this reason, we can obtain a better bound by replacing $m$ with $\min\{m,n+1\}$. 
    Since $\lp{(\x,z)}{2}\leq \sqrt{n+1}\normi{(\x,z)}\leq \sqrt{n+1}R$, part \ref{th:boundpart1} of the theorem follows.

    To show part \ref{th:boundpart2}, consider the formula $\Phi'(R):= \forall R'\in \R\{\Phi(R')\Longrightarrow (R'\leq R)\}$. Note that $\Phi'(R)$ can be written in prenex form as
    \begin{multline*}
    \forall R'\in \R \; \exists \La \in \R^m \; \forall (\t,w) \in \R^{n+1} \Bigg\{ \Bigg\{ \left[
    \lambda_1 \geq 0, \ldots, \lambda_m \geq 0, \quad \sum_{i=1}^m \lambda_i = 1
    \right]\wedge\\
     \Bigg[\left(\bigvee_{j=1}^n (t_j<0)\right) \vee (w<0) \vee 
    \left(\sum_{j=1}^n t_j \neq w\right) \bigvee_{j=1}^n \left[\left[
        \left( \sum_{i=1}^m \lambda_i f_{ij} \right)^r \neq w^s t_j^{r-s}
    \right]\wedge \left[
        -\left( \sum_{i=1}^m \lambda_i f_{ij} \right)^r \neq w^s t_j^{r-s}
    \right]\right] \\
    \vee \left(R' \left( \sum_{i=1}^m \lambda_i g_i - w \right) - \sum_{i=1}^m \lambda_i h_i > 0\right)
    \Bigg]\Bigg\} \vee (0\leq R' \leq R)\Bigg\}.
    \end{multline*}

    It is easy to see that $\Phi'(R)$ is satisfied if and only if 

    \begin{equation*}
        R\geq \max\{z : (\x,z) \text{ feasible for } {\bf F}_p\}.
    \end{equation*}

    Hence, we can apply Proposition \ref{prop:renegarbound} to $\Phi'(R)$ to conclude that, similarly to \eqref{eq:bound}, $\log R=\bits((m+n)r)^{O(mn)}.$ It remains to show that $m$ can be replaced by $\min\{m,n+1\}$. To this end, note that if the solution set of \ref{eq:OCFP} is bounded, then there exists a system $\f^T_i\x+g_iz\leq h_i, \; i\in I, \; \lp{\x}{p}\leq z$ with at most $n+1$ inequalities whose solution is still bounded. This is because the solution of problem \ref{eq:OCFP} is bounded if and only if the \emph{recession cone} of \ref{eq:OCFP} is trivial, namely,

    \begin{equation}\label{eq:recc}
        \K^n_p {\cap} \bigcap_{i=1}^m \H_i=\{\0\},
    \end{equation}

    where $\H_i$ is the halfspace $\{(\x,z)\in \R^{n+1} : \f^T_i\x+g_iz \leq 0\}$. 
    With $\Omega_1=\{(\x,z)\in \K^n_p : z=1\}$, \eqref{eq:recc} is equivalent to the emptiness of the intersection of the $(m+1)$ convex sets $\Omega_1, \H_1,\ldots, \H_m\subset \R^{n+1}$. 
    By Helly's theorem {(see, e.g. \cite{helly1923mengen})} there exists at most $(n+2)$ sets among $\Omega_1, \H_1,\ldots, \H_m\subset \R^{n+1}$ whose intersection is still empty. 
    The claim follows since $\Omega_1$ must be one of these sets.
    \end{proof}

    If the conjugate $q$ of $p$ is an even integer, for instance when $q=p=2$, the bound of Theorem \ref{th:bounds} can be improved as stated in the following result.
    \begin{theorem}\label{th:bounds2} 
    Let $r\in 2\N^*$ and $p=\frac{r}{r-1}$. Then
    \begin{enumerate}
        \item Any solution $(\x,z)$ of \ref{eq:OCFP}  satisfies $\norme{(\x,z)}\leq R$, where $\log R=\bits(r\min\{m,n\})^{O(\min\{m,n\})}$.
        \item Moreover, if the feasible set of \ref{eq:OCFP} is bounded, then the above bound holds for any feasible solution of \ref{eq:OCFP}.
    \end{enumerate}
\end{theorem}

    \begin{proof}
        First let us notice that $q= (1 - \frac{1}{p})^{-1} = (\frac{r}{r} - \frac{r-1}{r})^{-1} = r\in 2\N^*$. 
        The idea is then to adapt the proof of Theorem \ref{th:bounds} by considering the formula $\Phi(R)$ in the following standard form \ref{eq:sf}:
        \begin{multline*}
        \forall \La \in \R^m \; \exists w \in \R \Bigg\{ \Bigg\{ \left[
        \lambda_1 \geq 0, \ldots, \lambda_m \geq 0, \quad \sum_{i=1}^m \lambda_i = 1
        \right] \implies \\
        \Bigg[
        \sum_{j=1}^n \left( \sum_{i=1}^m \lambda_i f_{ij} \right)^r  = w^r,
        R \left( \sum_{i=1}^m \lambda_i g_i - w \right) - \sum_{i=1}^m \lambda_i h_i \leq 0, w\geq0
        \Bigg]\Bigg\} \wedge (R \geq 0)\Bigg\}.
        \end{multline*}
    
        One advantage given by the assumption on $p$ is that the variable $\t$ is not needed anymore, so $\Phi(R)$ is a standard formula \ref{eq:sf} of degree at most $r$ with $k=1$ free variable and $\omega=2$ quantifiers. Furthermore, $\Phi(R)$ consists of $O(m)$ atomic \rev{polynomial}  inequalities in $O(m)$ variables. The expansion of $\left( \sum_{i=1}^m \lambda_i f_{ij} \right)^r$ with integer coefficients of height at most $2^{\bits}$ yields coefficients with bit size at most $B=r(\bits+\log(m))$.
    
        From Proposition \ref{prop:renegarbound} it follows that $\Phi(R)$ can be satisfied by a positive number $R$ such that
    
        \begin{equation}
            \log R = r(\bits+\log(m))(O(m)r)^{O(m)}=\bits(mr)^{O(m)}.
        \end{equation}
    
        In the two parts of Theorem \ref{th:bounds} we can replace $m$ with $\min\{m,n+1\}$, yielding the desired result.
    \end{proof}


Let $R$ be the bound defined as

\begin{equation}\label{eq:R}
    R(p,n,m,\bits): =\begin{cases}
        \bits(r\min\{m,n\})^{O(n)} & \text{if } p=r\in 2\N^*,\\
        \bits(r\min\{m,n\})^{O(\min\{m,n\})} & \text{if } p=\frac{r}{r-1},\; r\in 2\N^*,\\
        \bits (rn)^{O(n\min\{m,n\})} & \text{otherwise,}
    \end{cases}
\end{equation}

and let $\K^n_p(R)=\{(\x,z)\in \K^n_p : z\leq R\}$. The \emph{discrepancy} of \ref{eq:OCFP} is the optimal value of the following optimization problem:

\begin{equation}\label{eq:discrepancy}
    \theta^*:=\min\{\theta : \f_i^T\x+g_iz\leq h_i+\theta, \; i\in M, \; (\x,z)\in \K^n_p(R)\}.
\end{equation}
Since $\K^n_p(R)$ is compact, the minimum in \eqref{eq:discrepancy} is always attained. 
In addition $\theta^*\leq 0$ if and only if \ref{eq:OCFP} is feasible.

The result below is the analog of ~\cite[Theorem 4.2]{porkolab1997complexity}.

\begin{theorem}\label{th:discrepancy}
Let $p=\frac{r}{s} \in \Q$, with $r > s\in \N^*$ and $\gcd(r,s)=1$. 
Assume that \ref{eq:OCFP} is infeasible. 
Then $-\log \theta^*=\bits(rn)^{O(n\min\{m,n\})}$.
If $r\in (2\Z)^*$ and $s=r-1$ then $-\log \theta^*=\bits(r\min\{m,n\})^{O(\min\{m,n\})}$.
\end{theorem}

The proof of Theorem \ref{th:discrepancy} is postponed at the end of Section~\ref{sec:generalocp}.


\section{Bounding the Complexity of the \texorpdfstring{$p$}{}-Order Feasibility Problem}\label{sec:comp}

{The result below is the analog of ~\cite[Lemma 5.2]{porkolab1997complexity}.}
\begin{theorem}\label{th:complexity}
If $p=\frac{r}{s} \in \Q$, with $r > s\in \N^*$ and $\gcd(r,s)=1$, then the feasibility of \ref{eq:OCFP} can be tested in $((m+n)r)^{O(n)}$ arithmetic operations over $\bits((m+n)r)^{O(n)}$-bit numbers.\\
If $p=r\in 2\N^*$ then the feasibility of \ref{eq:OCFP} can be tested in $(mr)^{O(n)}$ arithmetic operations over $\bits(mr)^{O(n)}$-bit numbers.\\
If $p=\frac{r}{r-1}$ with $r\in 2\N^*$ then the feasibility of \ref{eq:OCFP} can be tested in $(mr)^{O(m)}$ arithmetic operations over $\bits(mr)^{O(m)}$-bit numbers.
\end{theorem}

\begin{proof}
    By Lemma \ref{th:exreprs} \eqref{eq_gpcf_2}, the sentence

    \begin{multline}
        \exists (\x,z,\t)\in \R^{2n+1} \Bigg\{ \bigwedge_{i=1}^m (\f^T_i\x+g_iz\leq h_i),  t_1\geq0,\ldots, t_n\geq 0,
        \\
        \bigwedge_{j=1}^n \left[(x_j^r\leq z^{r-s}t_j^s)\wedge (-x_j^r\leq z^{r-s}t_j^s)\right], \dsum_{j=1}^nt_j\leq z\Bigg\}
    \end{multline}

    states that \ref{eq:OCFP} is feasible. From Theorem \ref{prop:renegarcomplexity} it follows that the validity of the above sentence can be determined in $((m+n)r)^{O(n)}$ operations over $\bits((m+n)r)^{O(n)}$-bit numbers. 

    If $p=r\in 2\N^*$ the formula becomes

    \begin{equation}
        \exists (\x,z)\in \R^{n+1}\left\{\bigwedge_{i=1}^m (\f^T_i\x+g_iz\leq h_i), \dsum_{j=1}^n x_j^r\leq z^r, z\geq 0\right\}.
    \end{equation}

    From Theorem \ref{prop:renegarcomplexity} it follows that the validity of the above sentence can be determined in $(mr)^{O(n)}$ operations over $\bits(mr)^{O(n)}$-bit numbers.

    Eventually, if $p=\frac{r}{r-1}$ with $r\in 2\N^*$, consider the sentence 

    \begin{equation} \label{eq:fease}
        \exists R\in \R \; \Phi(R),
    \end{equation}

    where $\Phi(R)$ is the formula defined in the proof of Theorem \ref{th:bounds2}. This sentence also states that \ref{eq:OCFP} is feasible. Observe that \eqref{eq:fease} consists of $O(m)$ polynomial inequalities of degree at most $r$ in $O(m)$ variables and has integer coefficients of bit size at most $B=r(\bits+\log(rm))$. The last part of the lemma follows again from Theorem \ref{prop:renegarcomplexity}.
\end{proof}

For $r=p=q=2$, Theorem \ref{th:complexity} directly implies the following corollary when considering the self-dual second order cone. 

\begin{corollary}\label{coro:socpcomplexity}
    The feasibility of an ${\rm SOCP}$ can be tested in $m^{O(\min\{m,n\})}$ arithmetic operations over $\bits m^{O(\min\{m,n\})}$.
\end{corollary}

As we will see in the next section, the bound of   Corollary \ref{coro:socpcomplexity} can be improved even further. 


\subsection{Complexity bounds on SOCP}

The result below is the analog of ~\cite[Theorem 5.1]{porkolab1997complexity}.
\begin{theorem}\label{th:socpcomplexity}
    The feasibility of an ${\rm SOCP}$ can be tested in $m\min\{m,n\}^{O(\min\{m,n\})}$ arithmetic operations over $\bits \min\{m,n\}^{O(\min\{m,n\})}$-bit numbers.
\end{theorem}

\begin{proof}
    If $m\leq n$ then the result follows directly from   Corollary \ref{coro:socpcomplexity}. 
    Let us assume that $m \geq n$. 
    
    Given a set $I\subset M=\{1,\ldots,m\}$, let us consider the following optimization problem
    \begin{equation}
        \theta(I):=\min\{\theta : \f^T_i\x+g_iz\leq h_i + \theta, \; i\in I, \; (\x,z)\in \K^n_2(R)\}, 
    \end{equation}
    where $R=R(2,n,m,\bits)$ \eqref{eq:R}. In particular, we have $\theta(M) = \theta^*$, the latter quantity has been defined in \eqref{eq:discrepancy}. 
    Denote by $(\x(I),z(I))$ the unique least norm solution of the system $\f^T_i\x+g_iz\leq h_i + \theta(I), i\in I, (\x,z)\in \K^n_2(R)$, and let  $V(I)=\{i\in I : \f^T_i\x(I)+g_iz(I)> h_i + \theta(I)\}$ \rev{be} the index set of constraints violated by $(\x(I),z(I))$. 
    A set $I$ is called a basis, if $V(J)\neq V(I)$ for any proper subset $J\subset I$. A basis $J$ is a basis for $I$, if $J\subseteq I$ and $V(J)=V(I)$. Any basis for $M$ is called optimal. In particular, if $S$ is an optimal basis, then

    \begin{equation}
        V(S)=V(M)=\emptyset, \text{ and consequently, } \theta(S)=\theta(M)=\theta^*.
    \end{equation}

    From Helly's theorem, it follows that $D:= \max\{|I| : I \text{ is a basis}\}\leq n+1$. Given an optimal basis $S$, we can apply Corollary \ref{coro:socpcomplexity} to $\f^T_i\x+g_iz\leq h_i, i\in S, (\x,z)\in \K^n_p$ and determine the feasibility of the original system that define the $\SOCP$ in $n^{O(n)}$ operations over $\bits n^{O(n)}$-bit numbers. Clarkson's algorithm~\cite{clarkson1995vegas} finds an optimal basis by performing expected $N=O(Dm+D^3\sqrt{m\log m}\log m)\leq m\rm{poly}(n)$ violation tests. Each of these checks whether $j\in V(I)$ for a sample set $I$ of cardinality $O(D^2\log D)$ and an index $j\in M\backslash I$. Note that the inclusion $j\in V(I)$ can be written as the sentence 

    \begin{multline}\label{eq:violated}
    V_I^j:= \forall (\x,z),(\y,w)\in \R^{n+1} \; \forall \theta,\rho \in \R \\\{ \{
    (\x^T\x\leq z^2, z\geq0)\wedge (\y^T\y\leq w^2, w\geq 0)\wedge S_I(\x,z,\theta)\\
     \wedge [S_I(\y,w,\rho)\Longrightarrow(\theta\leq \rho)]\wedge [S_I(\y,w,\theta)\Longrightarrow (\norme{(\x,z)}^2\leq \norme{(\y,w)}^2)]\}\\
    \Longrightarrow (\f^T_j\x+g_jz> h_j+\theta)\},
    \end{multline}

    where $S_I(\x,z,\theta)$ is the quantifier free formula $$\left\{\bigwedge_{i\in I} (\f^T_i\x+g_iz\leq h_i+\theta)\wedge (\norme{(\x,z)}^2\leq R^2)\right\}.$$ 

    Each violation test can thus be represented by a sentence in prenex form with $O(|I|)\leq O(n)$ polynomial inequalities of degree $2$ in $O(n)$ variables. Note also that the coefficient of these polynomial inequalities are integers of bit size $B\leq \max\{\bits,\log R\}=\bits\min\{m,n\}^{O(\min\{m,n\})}$. Now from Theorem \ref{prop:renegarcomplexity} it follows that each violation test can be accomplished in $n^{O(n)}$ operations over $\bits n^{O(n)}$-bit numbers. But the expected number of violation tests is bounded in $m\rm{poly}(n)$. Hence, we conclude that for all $n$ and $m$, testing the feasibility of a $\SOCP$ requires expected $mn^{O(n)}$ operations over $\bits n^{O(n)}$-bit numbers.
\end{proof} 

\begin{remark}
\label{rk:bpr}
Our complexity bound from Theorem~\ref{th:socpcomplexity} heavily relies on the result in  \cite{renegar1992computational} (Theorem \ref{prop:renegarcomplexity}). 
This latter result has been previously improved, e.g., in  \cite[Theorem~14.14]{basu2007algorithms}, where the authors obtain a similar algebraic complexity but a slight improvement in terms of bit size for the output and integers appearing in the intermediate computations. 
In particular the intermediate integers have bit sizes that are linear in the input bit size but does not depend on the number of polynomial (in)equalities. 
However in our case this does not yield any improvements because $\log R$ still depends polynomially on the number of (in)equalities. 
\end{remark}


In what follows we analyze the complexity bounds for other special cases for the values of $p$ in \ref{eq:OCFP}.

\begin{theorem}\label{th:complexeven}
    The feasibility of \ref{eq:OCFP} with $p=r\in 2\N^*$ can be tested in $m(r\min\{m,n\})^{O(n)}$ arithmetic operations over $\bits (r\min\{m,n\})^{O(n)}$-bit numbers.
\end{theorem}

\begin{proof}
    It is enough to follow the proof of Theorem \ref{th:socpcomplexity}, while relying on Corollary \ref{coro:socpcomplexity} instead of Theorem  \ref{th:complexity}, selecting $R=R(p,n,m,\bits)$ \eqref{eq:R}, and replacing the sentence \eqref{eq:violated} by the following one
    \begin{multline*}
    V_I^j:= \forall (\x,z),(\y,w)\in \R^{n+1} \; \forall \theta,\rho \in \R \\ \Bigg\{ \Bigg\{
    \left(\dsum_{k=1}^nx_k^r\leq z^r, z\geq0\right)\wedge \left(\dsum_{k=1}^ny_k^r\leq w^r, w\geq 0\right)\wedge S_I(\x,z,\theta)\\
     \wedge [S_I(\y,w,\rho)\Longrightarrow(\theta\leq \rho)]\wedge [S_I(\y,w,\theta)\Longrightarrow (\norme{(\x,z)}^2\leq \norme{(\y,w)}^2)]\Bigg\}\\
    \Longrightarrow (\f^T_j\x+g_jz> h_j+\theta)\Bigg\}.
    \end{multline*}
\end{proof}

\if{
\begin{theorem}\label{th:complexdualeven}
    Given $r\in 2\N^*$, the feasibility of \ref{eq:OCFP} with $p=\frac{r}{r-1}$ can be tested in  $m(r\min\{m,n\})^{O(\min\{m,n^2\})}$ arithmetic operations over $\bits (r\min\{m,n\})^{O(\min\{m,n^2\})}$-bit numbers.
\end{theorem}

\begin{proof}
    It is enough to follow the proof of Theorem \ref{th:socpcomplexity} while relying on Corollary \ref{coro:socpcomplexity} instead of Theorem  \ref{th:complexity}, and replacing the sentence \eqref{eq:violated} by
    \begin{multline*}
    V_I^j:= \forall (\x,z),(\y,w)\in \R^{n+1} \; \forall \theta,\rho \in \R \{ \{
    (\lp{\x}{p}\leq z)\wedge (\lp{\y}{p}\leq w)\wedge S_I(\x,z,\theta)\\
     \wedge [S_I(\y,w,\rho)\Longrightarrow(\theta\leq \rho)]\wedge [S_I(\y,w,\theta)\Longrightarrow (\norme{(\x,z)}^2\leq \norme{(\y,w)}^2)]\}\\
    \Longrightarrow (\f^T_j\x+g_jz> h_j+\theta)\}.
    \end{multline*}

    By Lemma \ref{th:exreprs} the membership of $(\x,z)$ in $\K^n_p$ can be expressed by the formula \eqref{eq_gpcf_2}. Then $V_I^j$ is a sentence involving $2$ quantifiers, $O(n)$ polynomial inequalities of degree at most $r$ in $O(n^2)$ variables. The value of $B$ remains unchanged. 
\end{proof}
}\fi

\subsection{Complexity bounds on general \texorpdfstring{\ref{eq:OCFP}}{}}
\label{sec:generalocp}
Once it is shown the complexity of SOCP, the proof can be modified to derive the complexities of more general cases.

\begin{theorem}\label{th:complexgc}
Let $p=\frac{r}{s} \in \Q$, with $r > s\in \N^*$ and $\gcd(r,s)=1$. 
The feasibility of \ref{eq:OCFP} can be tested in $\min\{[r(m+n)]^{O(n)}, m(rn)^{O(n^2)}\}$ arithmetic operations over $\bits\min\{[r(m+n)]^{O(n)}, (rn)^{O(n^2)}\}$-bit numbers.
\end{theorem}

\begin{proof}
It is enough to follow the proof of Theorem \ref{th:socpcomplexity}, while relying on Corollary \ref{coro:socpcomplexity} instead of Theorem  \ref{th:complexity}, selecting $R=R(p,n,m,\bits)$ \eqref{eq:R}, and replacing the sentence \eqref{eq:violated} by
    \begin{multline*}
    V_I^j:= \forall (\x,z),(\y,w)\in \R^{n+1} \; \forall \theta,\rho \in \R \{ \{
    (\lp{\x}{p}\leq z)\wedge (\lp{\y}{p}\leq w)\wedge S_I(\x,z,\theta)\\
     \wedge [S_I(\y,w,\rho)\Longrightarrow(\theta\leq \rho)]\wedge [S_I(\y,w,\theta)\Longrightarrow (\norme{(\x,z)}^2\leq \norme{(\y,w)}^2)]\}\\
    \Longrightarrow (\f^T_j\x+g_jz> h_j+\theta)\}.
    \end{multline*}

    By Lemma \ref{th:exreprs} the membership of $(\x,z)$ in $\K^n_p$ can be expressed by the formula \eqref{eq_gpcf_2}. 
    Then $V_I^j$ is a sentence involving $2$ quantifiers, $O(n)$ polynomial inequalities of degree at most $r$ in $O(n^2)$ variables. 
    The value of $B$ remains unchanged. 
    So for fixed dimension, Clarkson's algorithm yields a better bound.
\end{proof}
%
The result below is the analog of ~\cite[Theorem 5.4]{porkolab1997complexity}.
\begin{corollary}
\label{cor:complexdualeven}
Given $r\in 2\N^*$, the feasibility of \ref{eq:OCFP} with $p=\frac{r}{r-1}$ can be tested in  $m(r\min\{m,n\})^{O(\min\{m,n^2\})}$ arithmetic operations over $\bits (r\min\{m,n\})^{O(\min\{m,n^2\})}$-bit numbers.
\end{corollary}
\begin{proof}
This is a direct consequence of Theorem~\ref{th:complexgc} and Theorem~\ref{th:complexity}. 
\end{proof}

\begin{theorem}\label{th:height}
    Given an optimal basis $S$ of \eqref{eq:discrepancy}, in $(rn)^{O(n\min\{m,n\})}$ operations over $\bits(rn)^{O(n\min\{m,n\})}$-bit numbers we can find a system of univariate polynomial inequalities with integer coefficients such that $\theta^*$ is the only real solution of the system. In particular, $\theta^*$ is a root of a nontrivial polynomial $p(\theta)\in \Z[\theta]$ such that $\log \height (p)=\bits(r n)^{O(n\min\{m,n\})}$.
\end{theorem}

\begin{proof}
    Assume without loss of generality that the given basis $S$ coincides with $M$. In particular, $m\leq n+1$. From Von Neumman's minimax theorem and \eqref{eq:min2}, it follows that for $R\geq 0$
    \begin{multline*}
        \theta^*=\dmin_{(\x,z)\in \K^n_p(R)}\dmax_{\La\in \Delta_m} \dsum_{i=1}^m \lambda_i(\f_i^T\x + g_iz - h_i) \\
        =\dmax_{\La\in \Delta_m}\dmin_{(\x,z)\in \K^n_p(R)} \left\{\left(\dsum_{i=1}^m\lambda_i\f_i^T\right)\x + \left(\dsum_{i=1}^m \lambda_i g_i\right)z - \dsum_{i=1}^m \lambda_i h_i\right\} \\
        =\dmax_{\La\in \Delta_m}\left\{ \min\left[0,R\left(\dsum_{i=1}^m \lambda_i g_i - \lp{\dsum_{i=1}^m\lambda_i\f_i}{q}\right)\right]  - \dsum_{i=1}^m \lambda_i h_i \right\}.
    \end{multline*}
    Consider the formula
    \begin{equation*}
        \Lambda(\theta):= \forall \La\in \Delta_m \left\{ \min\left[0,R\left(\dsum_{i=1}^m \lambda_i g_i - \lp{\dsum_{i=1}^m\lambda_i\f_i}{q}\right)\right]  - \dsum_{i=1}^m \lambda_i (h_i+\theta)\leq 0 \right\},
    \end{equation*}
    where $R=R(p,n,m,\bits)$ \eqref{eq:R}. This formula states that $\theta\geq \theta^*$, and it can be written as follows
    \begin{multline*}
    \forall \La \in \R^m \; \exists (\t,w) \in \R^{n+1} \Bigg\{ \left[
    \lambda_1 \geq 0, \ldots, \lambda_m \geq 0, \quad \sum_{i=1}^m \lambda_i = 1
    \right] \implies \\
    \Bigg\{\Bigg[t_1\geq 0,\ldots, t_n\geq0, \bigwedge_{j=1}^n \left[\left[
        \left( \sum_{i=1}^m \lambda_i f_{ij} \right)^r = w^s t_j^{r-s}
    \right]\vee \left[
        -\left( \sum_{i=1}^m \lambda_i f_{ij} \right)^r = w^s t_j^{r-s}
    \right]\right], \\ \sum_{j=1}^n t_j = w \Bigg]\wedge 
    \left[\left(\sum_{i=1}^m \lambda_i (h_i+\theta)\geq 0 \right)\vee \left(R \left( \sum_{i=1}^m \lambda_i g_i - w \right) - \sum_{i=1}^m \lambda_i (h_i+\theta) \leq 0\right)\right] 
    \Bigg\}\Bigg\}.
    \end{multline*}
Now $\theta^*$ is the only real solution of $\Lambda^*(\theta):= \forall \theta'\in \R \{\Lambda(\theta)\wedge [\Lambda(\theta')\Longrightarrow (\theta \leq \theta^*)]\}$. 
By consecutively applying Theorem \ref{prop:renegarfree} to $\Lambda(\theta)$ and $\Lambda^*(\theta)$, the latter formula can be transformed into a quantifier free formula $\Lambda^{**}(\theta)$. 
This requires $[r(m+n)]^{O(mn)}\leq (rn)^{O(n\min\{m,n\})}$ operations with $\max\{\bits,\log R\}[r(m+n)]^{O(mn)}\leq \bits(rn)^{O(n\min\{m,n\})}$-bit numbers. 
The formula $\Lambda^{**}(\theta)$ is composed of univariate polynomial relations $p(\theta)\lhd 0$, where $\lhd\in\{\leq,<,=,\neq,>,\geq\}$. Since $\theta^*$ is the only real solution of $\Lambda^{**}(\theta)$, this formula can be transformed into an equivalent system of polynomial inequalities, which must contain a polynomial $p$ such that $p(\theta^*)=0.$
\end{proof}
\begin{theorem}\label{th:height2}
Given $r\in 2\N^*$, let $p=\frac{r}{r-1}$. 
Then we can replace the bounds of Theorem \ref{th:height} by $(r\min\{m,n\})^{O(\min\{m,n\})}$ and $\bits(r\min\{m,n\})^{O(\min\{m,n\})}$.
\end{theorem}

\begin{proof}
It is enough to follow the proof of Theorem \ref{th:height} by relying on Theorem \ref{th:bounds2} instead of Theorem \ref{th:bounds}, and writing the formula $\Lambda(\theta)$ as 
    \begin{multline*}
        \forall \La \in \R^m \; \exists w \in \R \Bigg\{ \left[
        \lambda_1 \geq 0, \ldots, \lambda_m \geq 0, \quad \sum_{i=1}^m \lambda_i = 1
        \right] \implies 
        \Bigg\{\Bigg[
        \sum_{j=1}^n \left( \sum_{i=1}^m \lambda_i f_{ij} \right)^r  = w^r, w\geq 0 \Bigg] \wedge \\\left[\left(\sum_{i=1}^m \lambda_i (h_i+\theta)\geq 0 \right)\vee \left(R \left( \sum_{i=1}^m \lambda_i g_i - w \right) - \sum_{i=1}^m \lambda_i (h_i+\theta) \leq 0\right)\right] \Bigg\}\Bigg\}.
        \end{multline*}
\end{proof}
\begin{remark}
    The minimal polynomial of an algebraic number $\alpha$ is the primitive irreducible polynomial $p\in \Z[x]$ such that $p(\alpha)=0$ and the leading coefficient of $p$ is positive. 
    The $\height$ of $\alpha$ is the height of its minimal polynomial. 
    Theorem \ref{th:height2} and the well-known Mignotte's inequality recalled in Proposition \ref{prop:height} show that $\log \height (\theta^*) = \bits(r\min\{m,n\})^{O(\min\{m,n\})}$.
\end{remark}
Theorem~\ref{th:height} and Theorem~\ref{th:height2} immediately imply Theorem \ref{th:discrepancy}. 
\begin{proof}[Proof of Theorem \ref{th:discrepancy}]
Since the problem \ref{eq:OCFP} is infeasible, one has $\theta^*>0$. 
Theorem \ref{th:height} and Theorem \ref{th:height2} imply that the positive algebraic number $\theta^*$ is a root of a nontrivial polynomial $p \in\Z[x]$ with integer coefficients of bit sizes $\bits(rn)^{O(n\min\{m,n\})}$ and $\bits(r\min\{m,n\})^{O(\min\{m,n\})}$, respectively. 
Proposition \ref{prop:height} implies that $\theta^*\geq \frac{1}{1+\height(p)}$, yielding the desired result.
\end{proof}
\begin{remark}
    If $p=r\in 2\N^*$, then one can write the formula $\Lambda(\theta)$ from the proof of Theorem \ref{th:height} as
    \begin{equation}
        \forall (\x,z)\in \R^{n+1}\left\{\bigwedge_{i\in M} [\f_i^T\x+g_iz-(h_i+\theta)\leq 0], \dsum_{j=1}^bx_j^r\leq z^r, z\geq 0, \norme{(\x,z)}^2\leq R^2\right\},
    \end{equation}
with $\log R\leq \bits(\min\{m,n\}r)^{O(n)}$. The bounds of Theorems \ref{th:discrepancy} and \ref{th:height} become $(\min\{m,n\}r)^{O(n)}$ and $\bits(\min\{m,n\}r)^{O(n)}$.
\end{remark}

\section{Generalization to several \texorpdfstring{$p$}{}-order cone constraints}
\label{sec:multi}
{
In this section, we generalize our results to the case of more than one $p$-order cone constraints. 
Namely, given $d \in \N^*$ and $n_1, \dots, n_d \in \N^*$, we consider the norm  constraints $(\x_k,z_k) \in \K_p^{n_k},  k=1,\ldots,d$. 
With $n:=n_1+\cdots+n_d$, let $(\x,\z):=(\x_1;\ldots;\x_d;z_1;\ldots;z_d) \in \R^{n+d}$ be the full vector of variables. 
\rev{Given $p_k=\frac{r_k}{s_k} \in \Q$, with $r_k > s_k\in \N^*$ and $\gcd(r_k,s_k)=1$, for $k=1, \ldots, d$, we consider the  $(p_1, \ldots, p_d)$-order feasibility problem that consists of determining whether there exists a real vector $(\x,\z) \in \R^{n+d}$ such that
\begin{align}
\label{eq:multiocp}
     F_1 \x_1 + G_1z_1+\cdots + F_d\x_d+ G_dz_d \leq H, &\nonumber\\
     (\x_k,z_k) \in \K_{p_k}^{n_k}, \; k=1,\ldots,d, &
\end{align}}
where $F_k \in \Z^{m\times n_k}, G_k, H \in \Z^{m}$ for $m \in \N^*$. 
We assume that all tuples $(\x_k,z_k)$ are independent. 
The above problem involves $n+d$ variables and $m+d$ constraints. 
If this independence assumption does not hold then one can always introduce additional variables satisfying it and encode the dependencies via linear equality constraints. 

The result of Theorem~\ref{th:complexity} can be generalized as follows. 

\rev{\begin{theorem}
\label{th:multicomplexity}
Let us consider \eqref{eq:multiocp} and assume that all tuples $(\x_k,z_k)$ are independent. Then,
\begin{itemize}
    \item[-] If $p_k=\frac{r_k}{s_k} \in \Q$, with $r_k > s_k\in \N^*$ and $\gcd(r_k,s_k)=1$ for $k=1, \ldots, d$, then the feasibility of \eqref{eq:multiocp} can be tested in $[(m+n+d)\max\{r_1, \ldots, r_d\}]^{O(n+d)}$ arithmetic operations over $\bits[(m+n+d) \max\{r_1, \ldots, r_d\}]^{O(n+d)}$-bit numbers.
\item[-] If $p_1 = \cdots = p_d=:p$ with $p=r\in 2\N^*$ then the feasibility of \eqref{eq:multiocp} can be tested in $[(m+d)r]^{O(n+d)}$ arithmetic operations over $\bits[(m+d)r]^{O(n+d)}$-bit numbers.
\item[-] If $p_1 = \cdots = p_d =: p$ with  $p=\frac{r}{r-1}$ and  $r\in 2\N^*$ then the feasibility of \eqref{eq:multiocp} can be tested in $[(m+d)r]^{O(md^2)}$ arithmetic operations over $\bits[(m+d)r]^{O(md^2)}$-bit numbers.
\end{itemize}
\end{theorem}
\begin{proof}
The argument of the proof is the same as the single constraint case, Theorem~\ref{th:complexity}. We emphasize below the main differences with respect to the single constraint case. First, by Lemma \ref{th:exreprs} \eqref{eq_gpcf_2}, the sentence

    \begin{multline*}
        \exists (\x,\z,\t)\in \R^{2n+d} \Bigg\{ \bigwedge_{i=1}^m \left(\dsum_{k=1}^d \left(\f_{ik}^T\x_k + g_{ik}z_k \right)\leq h_i\right),  \bigwedge_{k=1}^d\bigwedge_{j=1}^{n_k} (t_{kj}\geq 0),
        \\
        \bigwedge_{k=1}^d\bigwedge_{j=1}^{n_k} \left[(x_{kj}^{r_k}\leq z_k^{r_k-s_k}t_{kj}^{s_k})\wedge (-x_{kj}^{r_k}\leq z_k^{r_k-s_k}t_{kj}^{s_k})\right], \bigwedge_{k=1}^d\left(\dsum_{j=1}^{n_k}t_{kj}\leq z_k\right)\Bigg\}
    \end{multline*}

    states that \eqref{eq:multiocp} is feasible. From Theorem \ref{prop:renegarcomplexity} it follows that the validity of the above sentence can be determined in $[(m+n+d)\max\{r_1, \ldots, r_d\}]^{O(n+d)}$ operations over $\bits[(m+n+d) \max\{r_1, \ldots, r_d\}]^{O(n+d)}$-bit numbers.

    If $p_1 = \ldots = p_d := p$ and $p=r\in 2\N^*$ the formula becomes

    \begin{equation*}
        \exists (\x,\z)\in \R^{n+d}\left\{\bigwedge_{i=1}^m \left(\dsum_{k=1}^d \left(\f_{ik}^T\x_k + g_{ik}z_k \right)\leq h_i\right), \bigwedge_{k=1}^d\left(\dsum_{j=1}^{n_k}x_{kj}^{r}\leq z_k^{r}\right), \bigwedge_{k=1}^d (z_k\geq 0)\right\}.
    \end{equation*}

    From Theorem \ref{prop:renegarcomplexity} it follows that the validity of the above sentence can be determined in $[(m+d)r]^{O(n+d)}$ arithmetic operations over $\bits[(m+d)r]^{O(n+d)}$-bit numbers.
    
    On the other hand, thanks to the independence assumption, Equation \eqref{eq:min} in Lemma \ref{lemma:ptheq} can be generalized as follows:
    \begin{equation*}
        \dmin_{\genfrac{}{}{0pt}{}{(\x_k,z_k)\in \K^{n_k}_{p_k}:}{k=1,\ldots,d} }\left\{\dsum_{k=1}^d (\f_k^T\x_k + g_kz_k) : z_1=R_1,\dots,z_d=R_d\right\}=\dsum_{k=1}^dR_k(g_k-\lp{\f_k}{q_k}),
    \end{equation*}
where $\|\cdot\|_{q_k}$ is the dual of the $\|\cdot\|_{p_k}$ norm, i.e., $q_k$ is such that $\frac{1}{p_k} + \frac{1}{q_k} = 1$ for $k=1, \ldots, d$. Lemma \ref{th:boundseven} becomes $\norme{(\x,\z)}\leq R$, where $\log R = \bits [(\min\{m,n\}+d) \max\{r_1, \ldots, r_d\}]^{O(n+d)}$, meanwhile in Theorem \ref{th:bounds}, $\Theta(R)$ is replaced by  
\begin{align*}
        \Theta(R_1,\ldots, R_d)& := \dmin_{
        \genfrac{}{}{0pt}{}{(\x_k,z_k)\in \Omega_{R_k}:}{k=1,\ldots,d}
        }\max\left\{\dsum_{k=1}^d \left(\f_{1k}^T\x_k + g_{1k}z_k\right)-h_1,\ldots, \dsum_{k=1}^d \left(\f_{mk}^T\x_k + g_{mk}z_k\right)-h_m\right\}, \end{align*}
        
       where $\Omega_{R_k} :=\{(\x,z)\in \K^{n_k}_{\rev{p_k}} : z=R_k\}.$ Thus, for any $R_1,\ldots, R_d\geq 0$ and considering $\Delta_m$ as in \eqref{def:simplex}, we have that
    \begin{multline*}
        \Theta(R_1,\ldots, R_d)=\dmin_{\genfrac{}{}{0pt}{}{(\x_k,z_k)\in \Omega_{R_k}:}{k=1,\ldots,d}}\dmax_{\La\in \Delta_m} \dsum_{i=1}^m \lambda_i\left(\dsum_{k=1}^d \left(\f_{ik}^T\x_k + g_{ik}z_k \right)-h_i\right) \\
        =\dmax_{\La\in \Delta_m}\dmin_{\genfrac{}{}{0pt}{}{(\x_k,z_k)\in \Omega_{R_k}:}{k=1,\ldots,d}} \left[\dsum_{k=1}^d \left( \left(\dsum_{i=1}^m\lambda_i\f_{ik}^T\right)\x_k + \left(\dsum_{i=1}^m \lambda_i g_{ik}\right)z_k \right) - \dsum_{i=1}^m \lambda_i h_i\right] \\
        =\dmax_{\La\in \Delta_m}\left[ \dsum_{k=1}^d R_k\left(\dsum_{i=1}^m \lambda_i g_{ik} - \lp{\dsum_{i=1}^m\lambda_i\f_{ik}}{q_k}\right)  - \dsum_{i=1}^m \lambda_i h_i \right].
    \end{multline*}
Then $\Phi(R)$ is replaced by  
    \begin{equation*}
        \Phi(R_1,\ldots, R_d) := \forall \La \in \Delta_m \left\{ \dsum_{k=1}^d R_k\left(\dsum_{i=1}^m \lambda_i g_{ik} - \lp{\dsum_{i=1}^m\lambda_i\f_{ik}}{q_k}\right)  - \dsum_{i=1}^m \lambda_i h_i\leq 0 \wedge (R_1\geq 0, \ldots, R_d\geq 0)\right\}.
    \end{equation*}
    Hence, by Lemma \ref{th:exreprs} \eqref{eq_gpcf_1}, $\Phi(R_1,\ldots, R_d)$ can be rewritten in the standard form \ref{eq:sf} as follows
    
    \begin{multline*}
    \forall \La \in \R^m \; \exists (\t_k,w_k)_{k=1}^d \in \R^{n+d} \Bigg\{ \Bigg\{ \left[
    \lambda_1 \geq 0, \ldots, \lambda_m \geq 0, \quad \sum_{i=1}^m \lambda_i = 1
    \right] \implies \\
    \Bigg[
    \sum_{j=1}^{n_1} t_{1j} = w_1,\ldots, \sum_{j=1}^{n_d} t_{dj} = w_d \bigwedge_{k=1}^d\bigwedge_{j=1}^{n_k} \left[\left[
        \left( \sum_{i=1}^m \lambda_i f_{ikj} \right)^{r_k} = w_k^{s_k} t_{kj}^{r_k-s_k}
    \right]\vee \left[
        -\left( \sum_{i=1}^m \lambda_i f_{ikj} \right)^{r_k} = w_k^{s_k} t_{kj}^{r_k-s_k}
    \right]\right], \\
    \dsum_{k=1}^d R_k\left(\dsum_{i=1}^m \lambda_i g_{ik} - w_k\right) - \dsum_{i=1}^m \lambda_i h_i\leq 0 \bigwedge_{k=1}^d\bigwedge_{j=1}^{n_k} t_{kj}\geq 0
    \Bigg]\Bigg\} \wedge (R_1 \geq 0,\ldots, R_d\geq 0)\Bigg\}.
    \end{multline*}
yielding $\norme{(\x,\z)}\leq R$, with $\log R= \bits[(n+d) \max\{r_1, \ldots r_d\}]^{O(\min\{m,n+d\}(n+d)d)}$. 

Meanwhile, as in Theorem \ref{th:bounds2}, $\Phi(R_1,\ldots, R_d)$ can be rewritten in the standard form
 \begin{multline*}
        \forall \La \in \R^m \; \exists \w \in \R^d \Bigg\{ \Bigg\{ \left[
        \lambda_1 \geq 0, \ldots, \lambda_m \geq 0, \quad \sum_{i=1}^m \lambda_i = 1
        \right] \implies \Bigg[
        \bigwedge_{k=1}^d \sum_{j=1}^n \left( \sum_{i=1}^m \lambda_i f_{ikj} \right)^{r_k}  = w_{k}^{r_k},\\
        \dsum_{k=1}^d R_k\left(\dsum_{i=1}^m \lambda_i g_{ik} - w_k\right) - \dsum_{i=1}^m \lambda_i h_i\leq 0, w_1,\ldots, w_d\geq0
        \Bigg]\Bigg\} \wedge (R_1 \geq 0, \ldots, R_d\geq 0)\Bigg\},
        \end{multline*}
 yielding $\norme{(\x,\z)}\leq R$, with $\log R= \bits[(\min\{m,n\}+d) \max\{r_1, \ldots, r_d\}]^{O(\min\{m,n+d\}d^2)}$. The latter formula along with the sentence $\exists (R_1,\ldots, R_d)\in \R^d \Phi(R_1,\ldots, R_d)$, and  Theorem \ref{prop:renegarcomplexity} provide the remaining complexity bounds in case $p_1 = \cdots = p_d = p$ with  $p=\frac{r}{r-1}$ and  $r\in 2\N^*$. In such a case, the feasibility of \eqref{eq:multiocp} can be tested in $[(m+d)r]^{O(md^2)}$ arithmetic operations over $\bits[(m+d)r]^{O(md^2)}$-bit numbers.
\end{proof}}

Similarly, we can extend the results of Theorem \ref{th:socpcomplexity}, Theorem \ref{th:complexeven}, Corollary \ref{cor:complexdualeven} and Theorem \ref{th:complexgc}, respectively. 
\rev{\begin{theorem}
\label{th:multiothers}
Let us consider \eqref{eq:multiocp} and assume that all tuples $(\x_k,z_k)$ are independent. 
\begin{itemize}
    \item[-] If $p_1= \cdots = p_d=2$ (SOCP), the feasibility of \eqref{eq:multiocp} can be tested in $m[\min\{m,n\}+d]^{O(\min\{n+d,md^2\})}$ arithmetic operations over $\bits[\min\{m,n\}+d]^{O(\min\{n+d,md^2\})}$-bit numbers.
\item[-] If $p_1= \cdots = p_d=:r\in 2\N^*$, the feasibility of \eqref{eq:multiocp} can be tested in $m[r\min\{m,n\}+d]^{O(n+d)}$ arithmetic operations over $\bits [r\min\{m,n\}+d]^{O(n+d)}$-bit numbers.
\item[-] If $p_1 = \cdots = p_d =: p =\frac{r}{s}$ with $r\in 2\N^*$ and $p=\frac{r}{r-1}$, the feasibility of \eqref{eq:multiocp} can be tested in  either $[r(m+d)]^{O(md^2)}$ arithmetic operations over $\bits[r(m+d)]^{O(md^2)}$-bit numbers or $m[r(n+d)]^{O((n+d)n)}$ arithmetic operations over $\bits[r(n+d)]^{O(\min\{m,n+d\}d^2+ (n+d)n)}$-bit numbers.
\item[-] If $p_1 = \cdots =p_d =: p =\frac{r}{s} \in \Q$, with $r > s\in \N^*$ and $\gcd(r,s)=1$, the feasibility of \eqref{eq:multiocp} can be tested in 
either $[r(m+n+d)]^{O(n+d)}$ arithmetic operations over $\bits[r(m+n+d)]^{O(n+d)}$-bit numbers or $m[r(n+d)]^{O((n+d)n)}$ arithmetic operations over $\bits[r(n+d)]^{O((\min\{m,d\}+n)(n+d)d)}$-bit numbers. 
\end{itemize}
\end{theorem}}
}

\section{Applications}
\label{sec:applis}
The complexity of the feasibility problem studied in this paper has a direct impact in the complexity of solving, with some tolerance $\epsilon>0$, optimization problems involving $p$-order cones. For instance, if the objective function of the problem ranges in $[\texttt{lb}, \texttt{ub}]$, by binary search approaches, the complexity of solving the optimization problem equals $\log_2(\frac{\texttt{ub}-\texttt{lb}}{\epsilon})$ times the complexity of the feasibility oracle. In what follows, we provide detailed complexity results for some optimization problems of interest in different fields.
\subsection{Norm minimization}
We start with the classical application of $\ell_p$-norm minimization problem ($p$-NMP), that is defined in standard form as
\begin{align}\label{eq:NMP}\tag{$p$-NMP}
    \text{minimize  }  & \lp{\x}{p} \nonumber\\ 
    \text{subject to  } & A\x\leq \b, \nonumber
\end{align}
 where $A\in \Z^{m\times n}$ and $\b\in \Z^m$. Geometrically, this program stands for the vector closest to zero (\0) in the half-space of the normed space $(\R^n,\|\cdot \|_p)$ defined by $A\x\leq \b$. This problem has an extended representation to a $p$-order cone program by means of an auxiliary variable $z\in \R$ as below:
\begin{align}\label{eq:cp}
    \text{minimize  }  & z \\ 
    \text{subject to  } & A\x\leq \b, \nonumber\\
                        & (\x,z)\in \K^n_p.\nonumber
\end{align}
The feasibility problem associated to \eqref{eq:cp} corresponds to an instance of \ref{eq:OCFP}, with $\f_i=\a_i$ being the $i$-th column of $A$, $i=1,\ldots,m$;   $g_1=\cdots=g_m=0$, and $h_1=b_1,\ldots, h_m=b_m$. 
\begin{corollary}
\label{cor:normmin}
Let us assume that the coefficients of the input data $A,b$ have bit size at most $\bits$, and let $N \in \N, \varepsilon=2^{-N}$. 

If $p=\frac{r}{s}\in \Q$  with $r > s\in \N^*$ and $\gcd(r,s)=1$, the feasibility of \eqref{eq:cp} can be tested in $\min\{[r(m+n)]^{O(n)}, m(rn)^{O(n^2)}\}$ arithmetic operations over $\bits\min\{[r(m+n)]^{O(n)}, (rn)^{O(n^2)}\}$-bit numbers. 
If the problem is feasible, then according to the proof of Theorem \ref{th:bounds}, there exists a solution $(\x,z)$ satisfying $z\leq R$, with $\log R\leq \bits (rn)^{O(n\min\{m,n\})}$, yielding an upper bound for \ref{eq:NMP}.\\
An $\varepsilon$-optimal solution of \ref{eq:NMP} can be obtained through binary search in $ (\tau + N)\min\{[r(m+n)]^{O(n)}, m(rn)^{O(n^2)}\}$ arithmetic operations. 

For the euclidean norm minimization $2$-NMP, feasibility can be tested in $m\min\{m,n\}^{O(\min\{m,n\})}$ arithmetic operations over $\bits\min\{m,n\}^{O(\min\{m,n\})}$-bit numbers, and we obtain an upper bound of bit size $\bits\min\{m,n\}^{O(\min\{m,n\})}$ for $2$-NMP. \\
An $\varepsilon$-optimal solution of $2$-NMP can be obtained through binary search in $ (\tau + N) m\min\{m,n\}^{O(\min\{m,n\})}$ arithmetic operations.
\end{corollary}
\begin{proof}
The two feasibility estimates follow directly from Theorem \ref{th:complexgc} and Theorem \ref{th:socpcomplexity}, respectively. 
As mentioned at the beginning of this section, the optimization cost of the binary search procedure depends on the available lower and upper bounds on the minimum. 
In the case of norm minimization, we can obviously select 0 as a lower bound, and an upper bound bit size according to the one involved in the proof of Theorem \ref{th:bounds}. 
\end{proof}
Our complexity results can be similarly applied to minimize sum or maximum of norms, see \S~2.2 from \cite{alizadeh2003second} for the corresponding formulations as SOCPs in the case of Euclidean norms. 
{
\subsection{Support Vector Machines}\label{sec:svm}
A slightly different version of the norm minimization problem appears in supervised classification problems, in the so-called $\ell_p$-Support Vector Machines~\citep[see e.g.][]{blanco2020lp}. Given a training sample $\{(\x_1, y_1), \ldots, (\x_m, y_m)\} \subseteq \Q^m \times \{-1,+1\}$, the goal is to construct an hyperplane-based classifier separating the two classes ($-1$ and $1$) by maximizing the $\ell_p$-norm separation between them. The problem is stated as:
\begin{align}\label{prob:svm}\tag{$\ell_p$-SVM}
\text{minimize  } & \;\;\lp{{\bm \omega}}{p} + C \sum_{i=1}^m \xi_i\\
\text{subject to } & \; y_i({\bm \omega}^T x_i + b) \geq 1 - \xi_i, \quad \forall i=1, \ldots, m,\nonumber\\
& \; {\bm \omega} \in \R^n, b \in \R,\nonumber\\
& \; \xi \geq 0.\nonumber
\end{align}
where ${\bm \omega}$ and $b$ are the coefficients of the separating hyperplane ($\mathcal{H} = \{{\bm \omega}^T z + b: z \in \R^n\}$), and $\xi$ are the \rev{misclassification} errors. The parameter $C > 0$ allows one to find a trade-off between the margin separation and the \rev{misclassification}. 

Note that this problem can be rewritten in the shape of the above norm minimization \rev{problem} as follows:
\begin{align*}
\label{eq:svm}
\text{minimize  } & \;\; z + C \sum_{i=1}^m \xi_i \nonumber\\
\text{subject to } & \; y_i({\bm \omega}^T x_i + b) + \xi_i\geq 1 , \quad \forall i=1, \ldots, m,\nonumber\\
& \; {\bm \omega} \in \R^n, b \in \R,\\
& \; {\bm \xi} \geq 0,\nonumber\\
& \; ({\bm \omega}, z) \in \K_p^n\nonumber
\end{align*}
The complexity of testing feasibility of the above problem is very similarly to the case of norm minimization, stated in Corollary \ref{cor:normmin}. 
\begin{corollary}
\label{cor:svm}
Assume that the coefficients of the input data $(C, \x_i)$ are rational numbers with bit size at most $\bits$, and let $N \in \N, \varepsilon=2^{-N}$. \\
If $p=\frac{r}{s}\in \Q$  with $r > s\in \N^*$ and $\gcd(r,s)=1$, then the feasibility of \ref{prob:svm} can be tested in $\min\{[r(m+n)]^{O(n)}, m(rn)^{O(n^2)}\}$ arithmetic operations over $\bits\min\{[r(m+n)]^{O(n)}, (rn)^{O(n^2)}\}$-bit numbers. 
An $\varepsilon$-optimal solution of \ref{prob:svm} can be obtained through binary search in $ (\tau + N)\min\{[r(m+n)]^{O(n)}, m(rn)^{O(n^2)}\}$ arithmetic operations. \\
If $p=2$, it can be tested in $m\min\{m,n\}^{O(\min\{m,n\})}$ arithmetic operations over $\bits\min\{m,n\}^{O(\min\{m,n\})}$-bit numbers. \\
An $\varepsilon$-optimal solution of \ref{prob:svm} can be obtained through binary search in $ (\tau + N)m\min\{m,n\}^{O(\min\{m,n\})}$ arithmetic operations.
\end{corollary}
\subsection{Robust Least Squares}\label{sec_rls}
Least squares problems consist of finding the coefficient of a linear hyperplane that minimize the sum of the squares differences between the \emph{predicted}  and the \emph{observed} values. Then, given a dataset $\{(\x_1, y_1), \ldots, (\x_k, y_k)\} \in \R^n \times \R$, with input data $X$ and response data $\y$, a least square problem can be formulated as:
$$
\min_{{\bm \omega} \in \R^n} \|X {\bm \omega} - \y\|_2^2
$$
In case the data are uncertain, in Robust Least Squares~\citep{el1996robust,bertsimas2018characterization} allows to derive solutions to the system by incorporating uncertainty sets for the parameters $X$ and $\y$ in the above problem. Specifically, assuming that 
\begin{align*}
\mathcal{U}_X &= \{\widetilde{X} \in \R^{n\times k}: \|\widetilde{X}\|_p \leq \nu_X\}\\
\mathcal{U}_{\y} &= \{\widetilde{\y} \in \R^{k}: \|\widetilde{\y}\|_p \leq \nu_{\y}\}
\end{align*}
the Robust Least Squares method is stated as the following problem:
\begin{align*}
    \min_{\genfrac{}{}{0pt}{}{{\bm \omega} \in \R^n}{\widetilde{X} \in \mathcal{U}_X,  \widetilde{\y} \in \mathcal{U}_{\y}}}
    \|(X+\widetilde{X}) {\bm \omega} - (\y+\widetilde{\y})\|_2^2
\end{align*}
\cite{bertsimas2018characterization} proved that the above problem can be reformulated as follows:
\begin{align*}
    \min_{{\bm \omega} \in \R^n} \|X {\bm \omega} - \y\|_2^2 + \nu_X \|{\bm \omega}\|_q + \nu_{\y}
\end{align*}
where $\|\cdot\|_q$ is the dual of the $\|\cdot\|_p$ norm, i.e., $q$ is such that $\frac{1}{p} + \frac{1}{q} = 1$. Thus, the above problem can be reformulated as a $p$-OCP:
\begin{align*}
    \text{minimize } & \|X {\bm \omega} - \y\|_2^2 + \nu_X z + \nu_{\y}\\
    \text{subject to } &  {\bm \omega} \in \R^n,\\
    &z \in \R,\\
    &({\bm \omega},z) \in \K_q^n,
\end{align*}
which can be solved with a similar complexity as the one in Corollary \ref{cor:svm} from the previous subsection. 
\subsection{Continuous Locations Problems}\label{sec:weber}
One of the foundational problems in Facility Location is the Weber Problem~\citep{fekete2005continuous,weber1922ueber}. Given a set of points $\mathcal{X} = \{\x_1, \ldots, \x_d\} \subset \R^n$, the goal is to find a point $\x\in \R^n$ minimizing the sum of the ($\lambda$-weighted) distances to the points in $\mathcal{X}$. Using $\ell_p$-norms, the problem is stated as:
\begin{equation}\label{prob:weber}\tag{Weber}
\min_{\x\in \R^n} \sum_{i=1}^d \lambda_i \|\x-\x_i\|_p
\end{equation}
which can be equivalently rewritten as a standard conic $p$-OCP as follows:
\begin{align*}
    \text{minimize } & \sum_{i=1}^d \lambda_i z_i\\
    \text{subject to } & \w_i = \x_i - \x, \quad \forall i=1, \ldots, d, \\
    & \w_ i \in \R^n, \quad \forall i=1, \ldots, d,\\
    & \x \in \R^n,\\
    & (\w_i, z_i) \in \K_p^n, \quad \forall i=1, \ldots, d.
\end{align*}
A generalized version of the Weber problem is the Continuous Ordered Median Location Problem (COMP). Given weights $\omega_1, \ldots, \omega_d$ (one can assume without loss of generality that they are in $[-1,1]$), in the COMP, the distances from the points to the new points are sorted in non decreasing order, and the ${\bm \omega}$-weights are assigned to the sorted sequence of distances, i.e., the COMP can be formulated as:
\begin{equation}\label{prob:comp}\tag{COMP}
\min_{\x\in \R^n} \sum_{i=1}^d \omega_i \lambda_{(i)}\|\x-\x_{(i)}\|_p
\end{equation}
where $\x_{(i)} \in \{\x_1, \ldots, \x_{d}\}$ such that $\|\x-\x_{(1)}\| \geq \|\x-\x_{(2)}\| \geq \cdots \geq \|\x-\x_{(d)}\|$. 
This unified framework allows, by adequately choosing the ${\bm \omega}$-weights, to model different problem of interest, as constructing the point minimizing the maximum of the distances from $\mathcal{X}$ to the new point (${\bm \omega} = (1, 0, \ldots, 0)$, the sum of the $m$ largest distances (${\bm \omega}=(1, \stackrel{m}{\ldots}, 1, 0, \ldots, 0)$) and many other measures. In case $\omega_1 \geq \ldots \geq \omega_d \geq \omega_{d+1}:=0$, it is known~\citep{blanco2014revisiting} that this problem is convex and it can be rewritten as
\begin{align*}
    \text{minimize } & \sum_{i=1}^d  (u_i + v_i)\\
    \text{subject to } & \w_i = \x_i - \x, \quad \forall i=1, \ldots, d,\\
    & u_i + v_k \geq \omega_i \lambda_k z_k, \quad \forall i, k=1, \ldots, d,\\
    & \w_ i \in \R^n, \quad \forall i=1, \ldots, d,\\
    & \x \in \R^n,\\
    & \u, \v \in \R^n, \\
    & (\w_i, z_i) \in \K_p^n, \quad \forall i=1, \ldots, d. \\
\end{align*}
As a straightforward consequence of Theorem \ref{th:multiothers}, we obtain the following result. 
\begin{corollary}
\label{cor:weber}
Let $p=\frac{r}{s} \in \Q$, with $r > s\in \N^*$, $\gcd(r,s)=1$, and let $N \in \N, \varepsilon=2^{-N}$.
Assuming that the coefficients of the input data $(\x_i,\lambda_i)$ have bit size at most $\bits$, then the feasibility of \ref{prob:comp} can be tested in $(r(d^2+nd))^{O(nd)}$ arithmetic operations over $\bits(r(d^2+nd))^{O(nd)}$-bit numbers. 
An $\varepsilon$-optimal solution of \ref{prob:comp} can be obtained through binary search in $ (\tau + N)(r(d^2+nd))^{O(nd)}$ arithmetic operations. If $p=2$, it can be tested in $(nd)^{O(nd)}$ arithmetic operations over $\bits(nd)^{O(nd)}$-bit numbers, furthermore a  $\varepsilon$-optimal solution of \ref{prob:comp} can be obtained by binary search in $ (\tau + N)(nd)^{O(nd)}$ arithmetic operations.
\end{corollary}
Note that this complexity can be directly extended to mixed-norm continuous location problems where the distance to each of the points is measured with a different $\ell_p$-norm. 
In this case, we can obtain complexity estimates by considering the worse case scenario, as mentioned at the end of Section \ref{sec:multi}. Similar results can also be obtained for the multiple-allocation multiple-facility counterpart of the above problem that is described in \citep{blanco2016continuous}.
\subsection{Robust Linear Programming}
\label{sec:robustlp}
It has been shown in  \cite{ben1999robust} that the robust counterpart of a linear program with ellipsoidal uncertainties can be formulated as an SOCP. 
We use the notation from \S~3.2 in \cite{alizadeh2003second}. 
Let us consider the robust linear optimization problem:
\begin{align}
\label{eq:RLP}\tag{Robust-LP}
    \text{minimize  }  & \hat{\bc}^T  \hat{\x} \nonumber\\ 
    \text{subject to  } & \hat{A} \hat{\x} \leq \b, \nonumber \\
    & \x \in \R^n, \nonumber
\end{align}
where the constraint data $\hat{A} \in \R^{m \times n}$ and $\b \in \Z^m$ are not known exactly. 
To ease the presentation, the above problem can be rewritten as 
\begin{align}
\label{eq:RLP2}
    \text{minimize  }  & \bc^T \x \nonumber\\ 
    \text{subject to  } & A \x \leq \bm{0},\\
    & x_{n+1} = -1, \nonumber\\
    & \x \in \R^{n+1} \nonumber,
\end{align}
where $A=[\hat{A}, \b]$, $\bc = (\hat{\bc},0)$, $\x = (\hat{\x}, x_{n+1})$. 

In \cite{ben1999robust}, the authors consider the case where the uncertainty set is the Cartesian product of ellipsoidal regions, one for each row $\a_i^T$ of $A$  centered at some given row vector $\overline{\a}_i^T \in \Z^{n+1}$, namely in the set $\{\a_i \in \R^{n+1} : \a_i = \overline{\a}_i + B_i \bm{u}, \|\bm{u}\|_2 \leq 1  \}$, where $B_i$ is a positive semidefinite matrix. 
Then they show that the robust counterpart of \eqref{eq:RLP2} is 
\begin{align}
\label{eq:RLP3}
    \text{minimize  }  & \bc^T \x \nonumber\\ 
    \text{subject to  } & \|B_i \x \|_2 \leq - \overline{\a}_i^T \x, \quad i=1,\dots,m,\\
    & x_{n+1} = -1, \nonumber\\
    & \x \in \R^{n+1} \nonumber.
\end{align}
Thanks to Theorem \ref{th:multiothers}, the complexity of testing feasibility of \eqref{eq:RLP3} (resp.~optimizing) can be readily estimated, yielding the following result. 
\begin{corollary}
\label{cor:robust}
Assume that $B_i \in \Z^{(n+1) \times (n+1)}$ is positive definite for each $i=1,\dots,m$.
Let $\tau$ be the maximal bit size of the input data $(\bc, B_i^{-1},\overline{\a}_i)$, and let $N \in \N, \varepsilon=2^{-N}$. 
The feasibility of \ref{eq:RLP} can be tested in $(mn)^{O(mn)}$ arithmetic operations over $\bits (m n)^{O(m n)}$-bit numbers. 
An $\varepsilon$-optimal solution of \ref{eq:RLP} can be obtained through binary search in $ (\tau + N)(mn)^{O(mn)}$ arithmetic operations. 
\end{corollary}
\begin{proof}
For each $i=1,\dots,m$, let us introduce auxiliary variable $\y^{(i)} := B_i \x$, so that $\x = B_i^{-1} \y^{(i)}$, and $z_i:=- \overline{\a}_i^T B_i^{-1} \y^{(i)}$. 
Then each inequality constraint $\|B_i \x \|_2 \leq - \overline{\a}_i^T \x$ is equivalent to $\|\y^{(i)} \|_2 \leq z_i$. 
We have to consider the $(m-1)n$ linear equality constraints $B_i^{-1} \y^{(i)} = B_j^{-1} \y^{(j)}$, for all $i \neq j$. 
The corresponding SOCP involves $n_{\text{socp}} = n m$ variables $(\y^{(1)}, \dots, \y^{(m)})$ and $m_{\text{socp}} = m + (m-1)n$ linear equality constraints. 
The number $d$ of cone constraints is equal to $m$. 
By Theorem \ref{th:multiothers}, the feasibility of the resulting SOCP can be tested in $m_{\text{socp}}[\min\{m_{\text{socp}},n_{\text{socp}}\}+m]^{O(\min\{n_{\text{socp}}+m, m_{\text{socp}} m^2\})}$ arithmetic operations over $\bits[\min\{m_{\text{socp}},n_{\text{socp}}\}+m]^{O(\min\{n_{\text{socp}}+m,m_{\text{socp}} m^2\})}$-bit numbers.
The desired result follows after noticing that $\min\{m_{\text{socp}},n_{\text{socp}}\}+m = mn + m + \min\{0, m -n \} = O(mn)$ and $n_{\text{socp}} + m = m n + m \leq m_{\text{socp}} m^2$. 
\end{proof}
A similar complexity estimate can be obtained if some $B_i$ is not invertible.
In this case, for any $\x$ in the kernel of $B_i$, the corresponding cone constraint is replaced by a linear inequality constraint $0 \leq - \overline{\a}_i^T \x$. 
In the worse case scenario, the whole feasible set is the union of $2^{m}$ feasible regions obtained by splitting the ambient space into the range and kernel of the matrices $B_1,\dots,B_m$. 
For the sake of simplicity, we restrict ourselves to the case of ellipsoidal uncertainties but the extension to $\ell_p$-ball based uncertainties might be derived similarly using the robust reformulations in \cite{bertsimas2004robust}. 

\section{Conclusions}

In this paper, we have addressed fundamental questions related to the complexity of solving mathematical optimization problems involving $\ell_p$-norms. 
While these problems can be equivalently reformulated as second-order cone (SOC) problems and further as semidefinite programming (SDP) problems, we demonstrate that leveraging the explicit structure of $p$-order cones ---particularly SOC--- yields improved complexity bounds compared to following the full reformulation path to SDP. 
Furthermore, specific choices of $p$ allow one for even more refined complexity results. 
We also investigate upper bounds for the norm of a solution when the problem is known to be feasible, as well as analyze the discrepancy in cases of infeasibility. 
The implications of our findings are explored in applications such as $\ell_p$-Support Vector Machines, robust optimization with $\ell_p$-norm-based uncertainty sets, and single-facility ordered continuous location problems involving $\ell_p$-norms, all of which have seen significant recent interest.

This work contributes to the deeper understanding of $p$-order cones and their role in solving optimization problems that involve these structures. Many optimization problems incorporating $\ell_p$-norm constraints also include integer or binary decision variables to model on/off or disjunctive constraints. While these problems lose convexity, efficiently solving their continuous relaxations has been key to developing branch-and-bound algorithms for Mixed Integer Second Order Cone Optimization (MISOCO) problems, enabling the solution of reasonably sized instances. However, Mixed Integer $p$-Order Cone Optimization (MI$p$OCO) problems have so far been approached only through reformulations as MISOCO problems. Extending our understanding and methods to address MI$p$OCO problems directly will be a promising direction for future research.

\rev{Another interesting research direction stemming from this work concerns the sharpness of the provided upper bounds. While we have established upper bounds on the number of arithmetic operations and the bit complexity required to decide the feasibility of $p$-order cone constraints, it remains an open question whether these bounds are tight. To the best of our knowledge, there are no matching lower bounds in the literature for the general $p$-order feasibility problem. For specific values of $p$, such as $p=2$, related problems, such as second-order feasibility and semidefinite feasibility, are known to admit only weakly polynomial algorithms, and no strongly polynomial algorithms are known. These results suggest intrinsic computational limitations, but a formal characterization of lower bounds, particularly for general rational values of $p$, remains elusive. }

\section*{Acknowledgments}
This work benefited from the HORIZON–MSCA-2023-DN-JD of the European Commission under the Grant Agreement No 101120296 (TENORS), the AI Interdisciplinary Institute ANITI funding, through the French ``Investing for the Future PIA3'' program under the Grant agreement n${}^\circ$ ANR-19-PI3A-0004 as well as the National Research Foundation, Prime Minister’s Office, Singapore under its Campus for Research Excellence and Technological Enterprise (CREATE) programme. 
This research has also been partially supported by grant PID2020-114594GB-C21 funded by MICIU/AEI/ 10.13039/501100011033, grant RED2022-134149-T funded by MICIU/AEI /10.13039/501100011033(Thematic Network on Location Science and Related Problems), and the IMAG-María de Maeztu grant CEX2020-001105-M/AEI/10.13039/501100011033.

The research stay of the third author at LAAS CNRS was partly funded by the LabEx CIMI (ANR-11-LABX-0040). 

The authors thank the anonymous referees for their constructive comments and valuable suggestions, which helped improve the clarity and quality of previous versions of the manuscript.



\end{document}